\documentclass{article}

\usepackage{amsmath}
\usepackage{amssymb}
\usepackage{commath}

\usepackage{float}
\usepackage{multirow}

\usepackage{xcolor}

\usepackage{titlesec}
\usepackage[section]{placeins}

\usepackage[utf8]{inputenc}
\usepackage[english]{babel}

\usepackage{graphicx}
\usepackage[labelformat=simple]{subcaption}

\usepackage{indentfirst}
\usepackage{enumitem}

\usepackage{amsmath}
\usepackage{amssymb}
\usepackage{commath}

\usepackage[numbers,square]{natbib}
\usepackage{float}
\usepackage{hyperref}
\usepackage[all]{hypcap}
\usepackage{multirow}

\usepackage{xcolor}

\usepackage{titlesec}
\usepackage[section]{placeins}
\titleformat{\paragraph}
{\normalfont\normalsize\bfseries}{\theparagraph}{1em}{}
\titlespacing*{\paragraph}
{0pt}{3.25ex plus 1ex minus .2ex}{1.5ex plus .2ex}

\usepackage{etoolbox}

\newtheorem{definition}{Definition}

\newcommand{\me}{\mathrm{e}}  
\newcommand{\iu}{{i\mkern1mu}} 
\newcommand{\rpm}{\sbox0{$1$}\sbox2{$\scriptstyle\pm$} 
  \raise\dimexpr(\ht0-\ht2)/2\relax\box2 }
  
  \DeclareMathOperator{\sign}{sign}
\DeclareMathOperator{\arctantwo}{arctan2} 
  
\makeatletter 
\renewenvironment{cases}[1][l]{\matrix@check\cases\env@cases{#1}}{\endarray\right.}
\def\env@cases#1{%
  \let\@ifnextchar\new@ifnextchar
  \left\lbrace\def\arraystretch{1.2}%
  \array{@{}#1@{\quad}l@{}}}
\makeatother

\pretocmd{\eqref}{Eq.~}{}{}
\newcommand\defref{Definition~\ref}

\newcommand{\Secref}[1]{Section~\ref{#1}} 
\newcommand{\tableref}[1]{Table~\ref{#1}}

\newcommand{\R}{\mathbb{R}}

\usepackage[toc,page]{appendix}

\usepackage{authblk}

\title{Isoptics of log-aesthetic curves}

\author[1]{Ferenc Nagy}

\affil[1]{University of Debrecen\\ Faculty of Informatics\\ 26 Kassai St.\\ 4028
  Debrecen\\ Hungary}
\affil[1]{University of Debrecen\\ Doctoral School of Informatics\\ 26 Kassai St.\\ 4028
  Debrecen\\ Hungary}

\date{January 2021}

\begin{document}

\maketitle

\begin{abstract}
The log-aesthetic curve has a significant factor in the field of aesthetic design to meet the high industrial requirements. It has much exceptional property and a large number of research papers are published, since its introduction. It can be parametrized by tangential angle, and it is autoevolute, meaning that its evolute is also a log-aesthetic curve. In this paper, we will examine the isoptics of the log-aesthetic curve and determine if it is autoisoptic as well. We will show that the log-aesthetic curve only coincides its isoptic in some specific cases.
We will also provide a new formula to calculate the isoptic curve of any tangential angle parametrized curves.
\end{abstract}

\section{Introduction}

To design highly aesthetic objects, the use of high-quality curves is very
important. The main property of such curves is the curvature, which should be monotonically varying.
Among them, the log-aesthetic stated as the most promising for aesthetic design
 \cite{levien2009interpolating} because its curvature and torsion graphs are both linear \cite{yoshida2009log}.

The log-aesthetic curve is originated from Harada et al. \cite{harada1997study, harada1999aesthetic}. They revealed that natural aesthetic curves have such a property that their logarithmic distribution diagram of curvature (LDDC) can be approximated by straight lines. Based on their work Miura et al. \cite{miura2005derivation,miura2006general} have defined the logarithmic curvature graph (LCG), an analytical version of the LDDC as the following: when a curve (given with arc length $s$ and radius of curvature $\rho$) is subdivided into infinitesimal segments such that $\Delta\rho/\rho$ is constant, the LCG represents the relationship between $\rho$ and $\Delta s$ in a double logarithmic graph. The authors formulated the curve whose LCG is strictly expressed by a straight line (with slope $\alpha$). Such curves satisfy the following equation:
\begin{equation} \label{eq:fundamental}
	\log\rho\frac{\dif s}{\dif \rho} = \alpha \log \rho + c,
\end{equation}
where $c$ is a constant.
\eqref{eq:fundamental} is the fundamental equation of log-aesthetic curves.

\begin{figure}[htb]
\centering
  \includegraphics[width=0.9\linewidth]{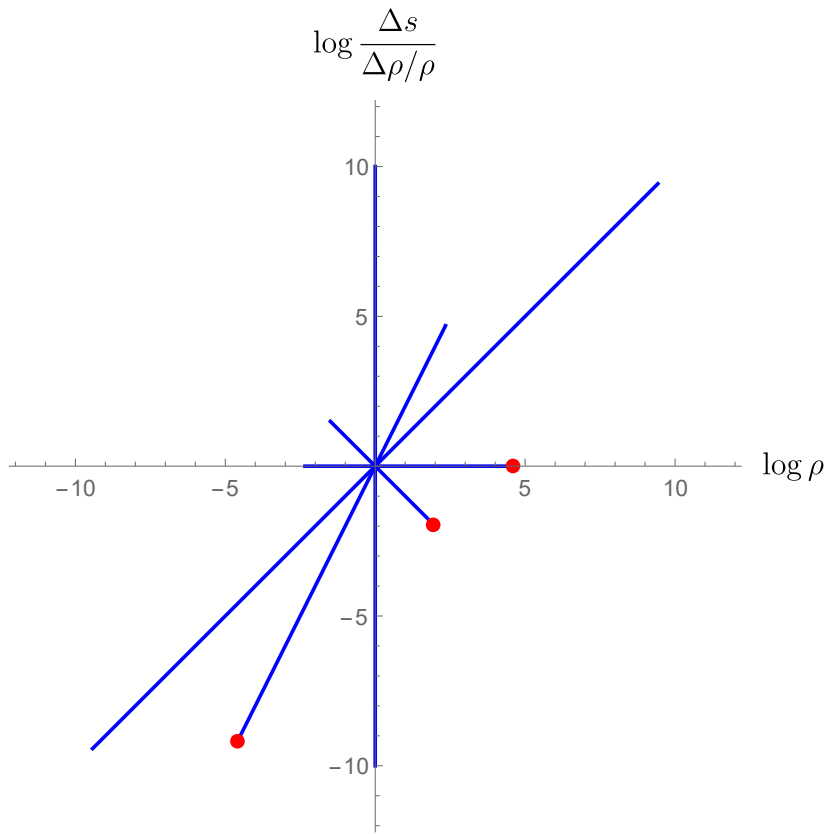}
\caption{
The logarithmic curvature graph of several log-aesthetic curves with different $\alpha$ slopes: clothoid ($\alpha=-1$), Nielsen's spiral ($\alpha=0$), logarithmic spiral ($\alpha=1$), circle involute ($\alpha=2$), circle ($\alpha=\rpm\infty$).
The curves are generated using parametric equations (on the same $[-3\pi,3\pi]$ interval), which parameter is limited  (red dots), except in the cases of $\alpha=1$ and $\alpha=\rpm\infty$.}
\label{fig:LCG_all}
\end{figure}

The log-aesthetic curve can also be considered as the generalization of the clothoid, logarithmic spiral, Nielsen's spiral, and the circle involute. The LCG of these curves is a straight line with different slopes that can be drawn as log-aesthetic curves with different shape parameter $\alpha$ (see \figref{fig:LCG_all}).
The general equations of the log-aesthetic curve will be presented in \Secref{sec:log-aesthetic_formulas}.

Besides the linear curvature, the log-aesthetic curve has another interesting property. It is autoevolute, that means the evolute of a log-aesthetic curve with shape parameter $\alpha$ is also a log-aesthetic curve with shape parameter $-1/(\alpha-2)$ \cite{yoshida2012evolutes}.

This characteristic led us to study the isoptic curve of the log-aesthetic curve.
The isoptic curve of a curve is constructed by the involving lines with
a given angle intersect each other at a certain point of the isoptic curve. Some results of calculating the isoptics of various plane curves will be presented in \Secref{sec:isoptic_curve}.
The similar attribute to autoevolute is called the autoisoptic, meaning that the given curve coincide its isoptic curve. In \cite{boris2019examples}, there are some autoisoptic curves, among the logarithmic spiral, which is the log-aesthetic curve with $\alpha=1$. It is also evident that the isoptic of a circle is also a circle (log-aesthetic curve with $\alpha=\pm\infty$).

In \Secref{sec:isoptic_log-aesthetic}, we will present a new formula for calculating the isoptic curve of the log-aesthetic curve (and for all the tangential angle parametrized curves). We will also prove that for arbitrary $\alpha \in \mathbb{R}$ the log-aesthetic curve is not autoisoptic.

\section{Previous results}

In this section, we are going to summarize the work of Yoshida and Saito \cite{yoshida2006interactive}. They analyzed the properties of the log-aesthetic curve and derived a general formula to draw it. Their related results will be presented in the first subsection. Besides, in the second part, we will introduce the isoptic curves and show some results in the topic.

\subsection{Generating log-aesthetic curves} \label{sec:log-aesthetic_formulas}

The general formula is obtained from the relationship between the arc length and the radius of curvature of the log-aesthetic curve, using a reference point (at the origin) and placing the following constraints on it: the tangent vector is $\big[1$ $0\big]^T$, and the radius of curvature is $1$. By transforming the curve such that the above constraints are satisfied, a point $P(\theta)$ of the log-aesthetic curve whose tangential angle is $\theta$ can be expressed on the complex plane as (\cite{yoshida2006interactive}):
\begin{equation} \label{eq:log-aesthetic_point_theta}
    P(\theta)= 
\begin{cases}
    \int_{0}^{\theta}  \me^{(1+\iu)\Lambda\psi}\  d\psi
        & \text{if } \alpha = 1\\
    \int_{0}^{\theta}  \big((\alpha-1)\Lambda\psi+1\big)^{\frac{1}{\alpha-1}} \me^{\iu\psi}\  d\psi
        & \text{otherwise},
\end{cases}
\end{equation}
where $\alpha \in \R$ and $\Lambda \in \R^{+}$ are parameters. $\alpha$ is the slope of the LCG.
When $\alpha \neq 1$, all the aesthetic curves are congruent under similarity transformations without depending on the value of $\Lambda(\neq 0$).
Since similarity transformations preserve angles, in our examples we focus on the cases of $\Lambda=1$.

The radius of curvature as the function of $\theta$ can be expressed as (\cite{yoshida2006interactive}):
\begin{equation} \label{eq:radius_curvature_theta}
    \rho(\theta)= 
\begin{cases}
    \me^{\Lambda\theta}
        & \text{if } \alpha = 1\\
    \big((\alpha-1)\Lambda\theta+1\big)^{\frac{1}{\alpha-1}} 
        & \text{otherwise}.
\end{cases}
\end{equation}
When $\theta=0$, $\rho=1$. The function increases monotonically when $\Lambda \neq 0$. In case of $\Lambda = 0$, $\rho$ is constant $1$ and the log-aesthetic curve is a circle. Similarly, when $\alpha = \rpm \infty$ because taking the limit of \eqref{eq:radius_curvature_theta} when $\alpha$ approaches $\rpm \infty$ we get $\rho=1$.

The tangential angle $\theta$ and the arc length $s$ are related by (\cite{yoshida2006interactive}):
\begin{equation} \label{eq:theta_arc}
    \theta(s) = 
\begin{cases}
    \frac{1-\me^{-\Lambda s}}{\Lambda}  
        & \text{if } \alpha = 0\\
    \frac{\log(\Lambda s + 1 )}{\Lambda}
        & \text{if } \alpha = 1\\
    \frac{(\Lambda\alpha s +1)^{(1-\frac{1}{\alpha})}-1} {\Lambda(\alpha-1)}  
        & \text{otherwise}.
\end{cases}
\end{equation}
Therefore, a point $C(s)$ on the aesthetic curve whose arc length is $s$, can be defined on the complex plane as (\cite{yoshida2006interactive}):
\begin{equation} \label{eq:log-aesthetic_point_arc}
    C(s)= 
\begin{cases}
    \int_{0}^{s}\exp({\iu\frac{1-\me^{-\Lambda u}}{\Lambda}})\ d u
        & \text{if } \alpha = 0\\
    \int_{0}^{s}\exp({\iu\frac{\log(\Lambda u + 1 )}{\Lambda}})\ d u
        & \text{if } \alpha = 1\\
    \int_{0}^{s}\exp({\iu\frac{(\Lambda\alpha u +1)^{(1-\frac{1}{\alpha})}-1} {\Lambda(\alpha-1)}})\  d u
        & \text{otherwise}.
\end{cases}
\end{equation}
\eqref{eq:log-aesthetic_point_arc} and \eqref{eq:log-aesthetic_point_theta} represent the same curve.
Using \eqref{eq:theta_arc}, the radius of curvature can also be expressed from the arc length $s$ (\cite{yoshida2006interactive}):
\begin{equation} \label{eq:radius_curvature_arc}
    \rho(s)= 
\begin{cases}
    \me^{\Lambda s}
        & \text{if } \alpha = 0\\
    \big(\Lambda\alpha s + 1 \big)^{\frac{1}{\alpha}} 
        & \text{otherwise.}
\end{cases}
\end{equation}

Since $\rho$ can change from $-\infty$ to $+\infty$, the tangential angle $\theta$ and arc length $s$ may have upper or lower bound depending on the value of $\alpha$ because of the negative bases of the fractional exponents (see \tableref{table:arc_theta_bounds}).

\begin{figure}[h!]
  \captionsetup[subfigure]{justification=centering}
	\centering
  \begin{subfigure}[t]{0.22\textwidth}
		\centering
		\includegraphics[width=\textwidth]{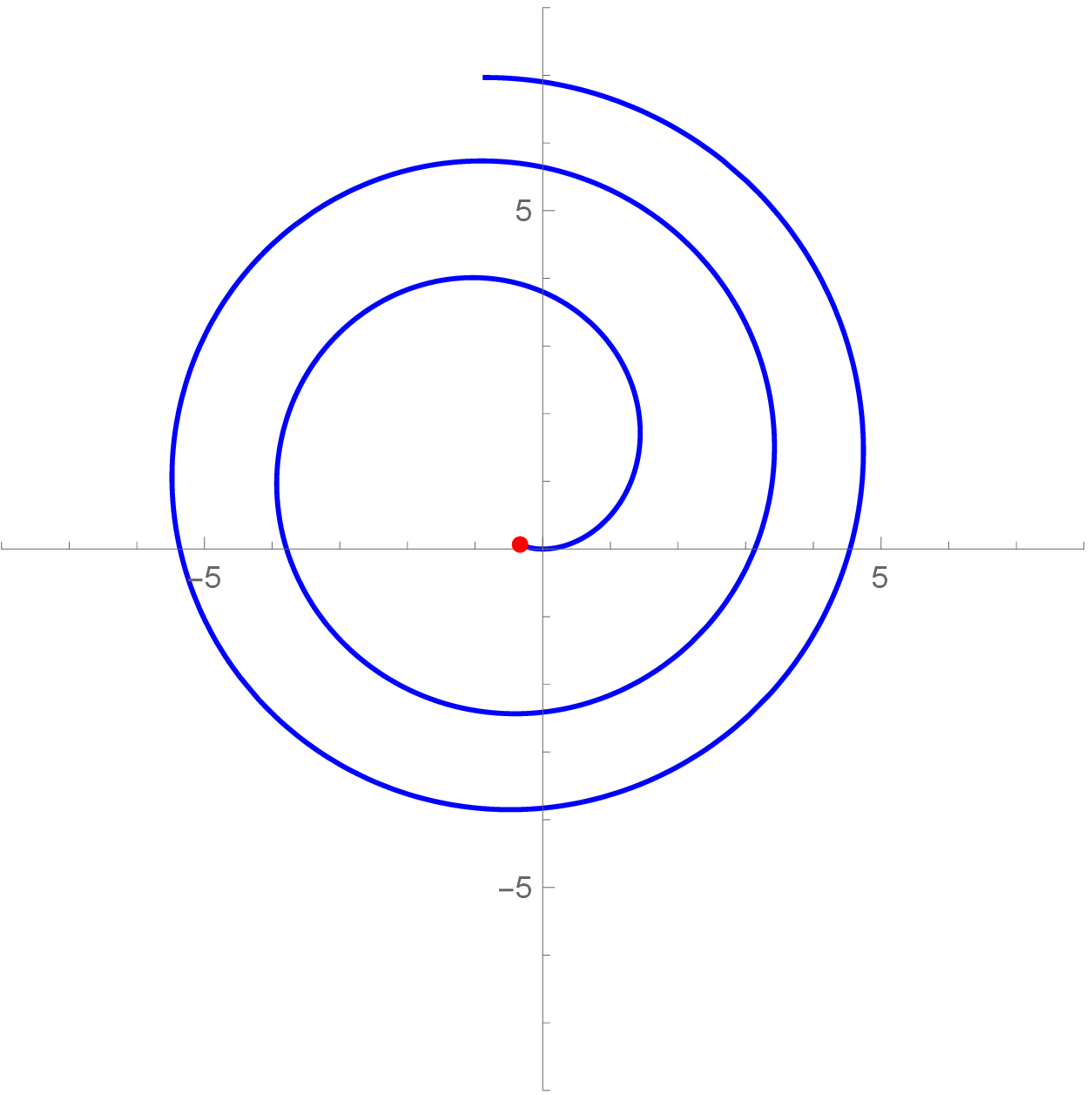}
		\caption{$\alpha=3$}\label{fig:log-aesthetic_alpha_3}
	\end{subfigure}
 \hfill
\begin{subfigure}[t]{0.22\textwidth}
		\centering
		\includegraphics[width=\textwidth]{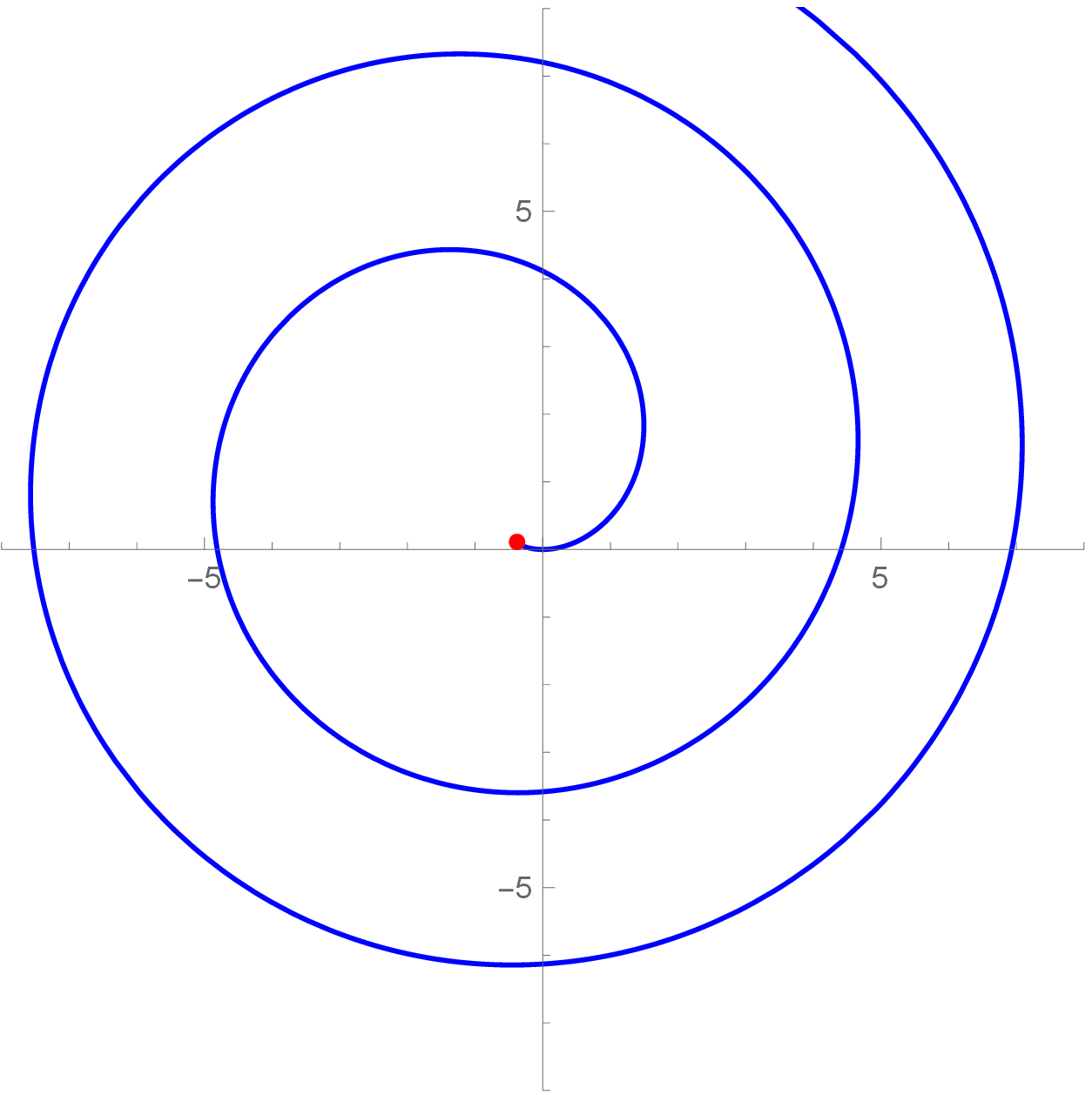}
		\caption{$\alpha=2.5$}
	\end{subfigure}
\hfill
\begin{subfigure}[t]{0.22\textwidth}
		\centering
		\includegraphics[width=\textwidth]{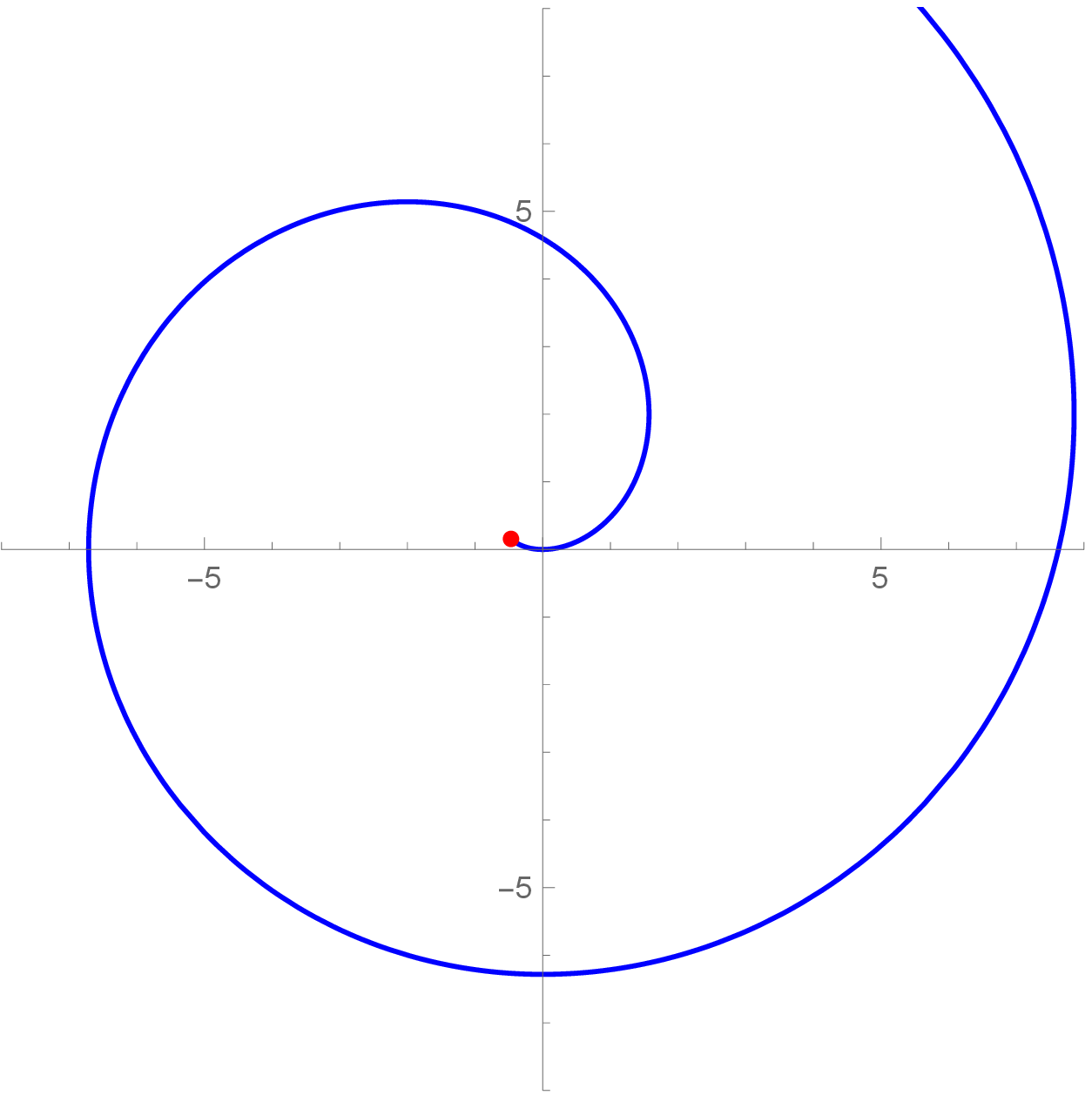}
		\caption{$\alpha=2$,\\Circle involute}
	\end{subfigure}
 \hfill
	\begin{subfigure}[t]{0.22\textwidth}
		\centering
		\includegraphics[width=\textwidth]{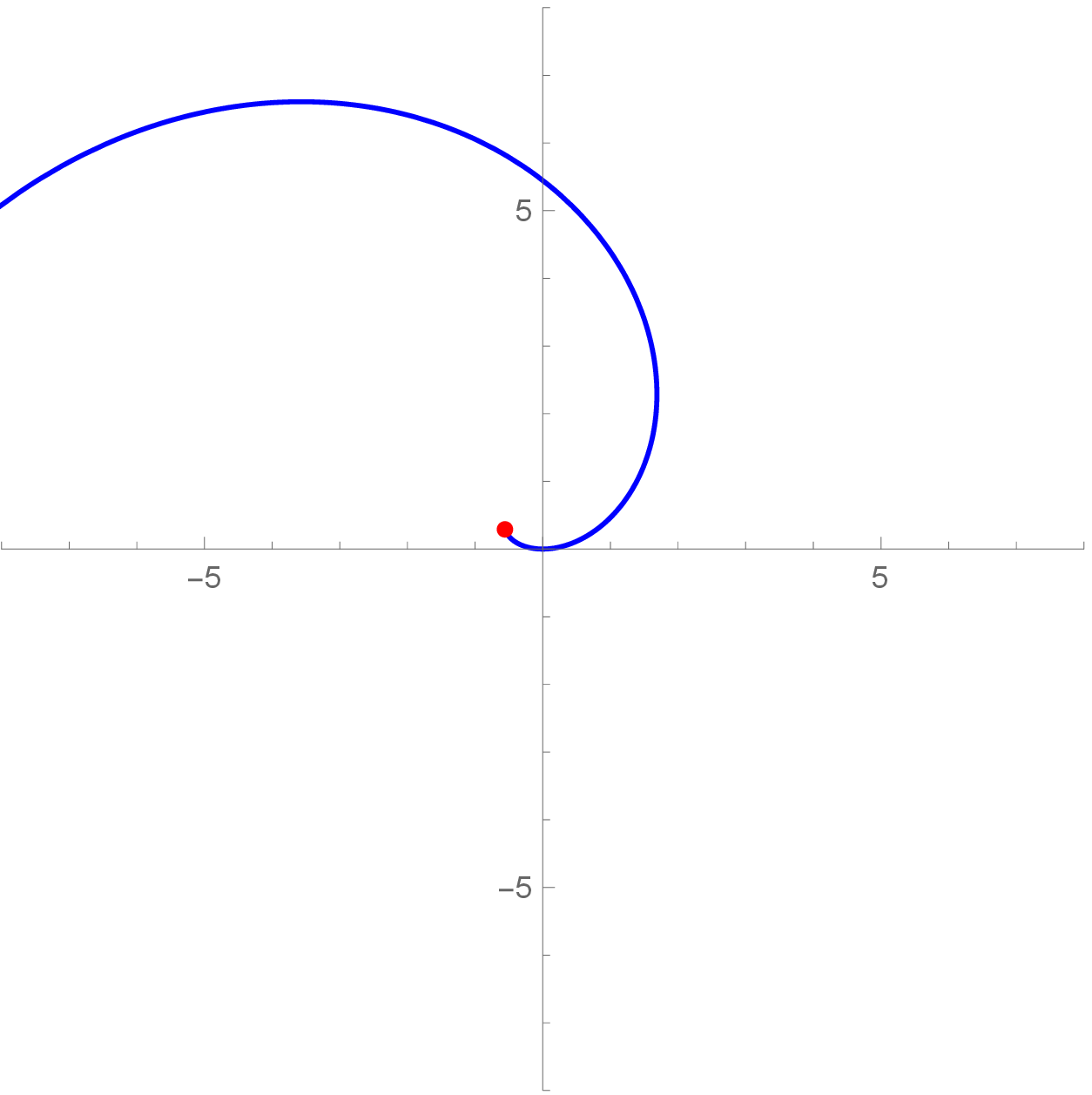}
		\caption{$\alpha=1.5$}\label{fig:log-aesthetic_alpha_1_5}
	\end{subfigure}
  \\[2em]
	\centering
	\begin{subfigure}[t]{0.22\textwidth}
		\centering
		\includegraphics[width=\textwidth]{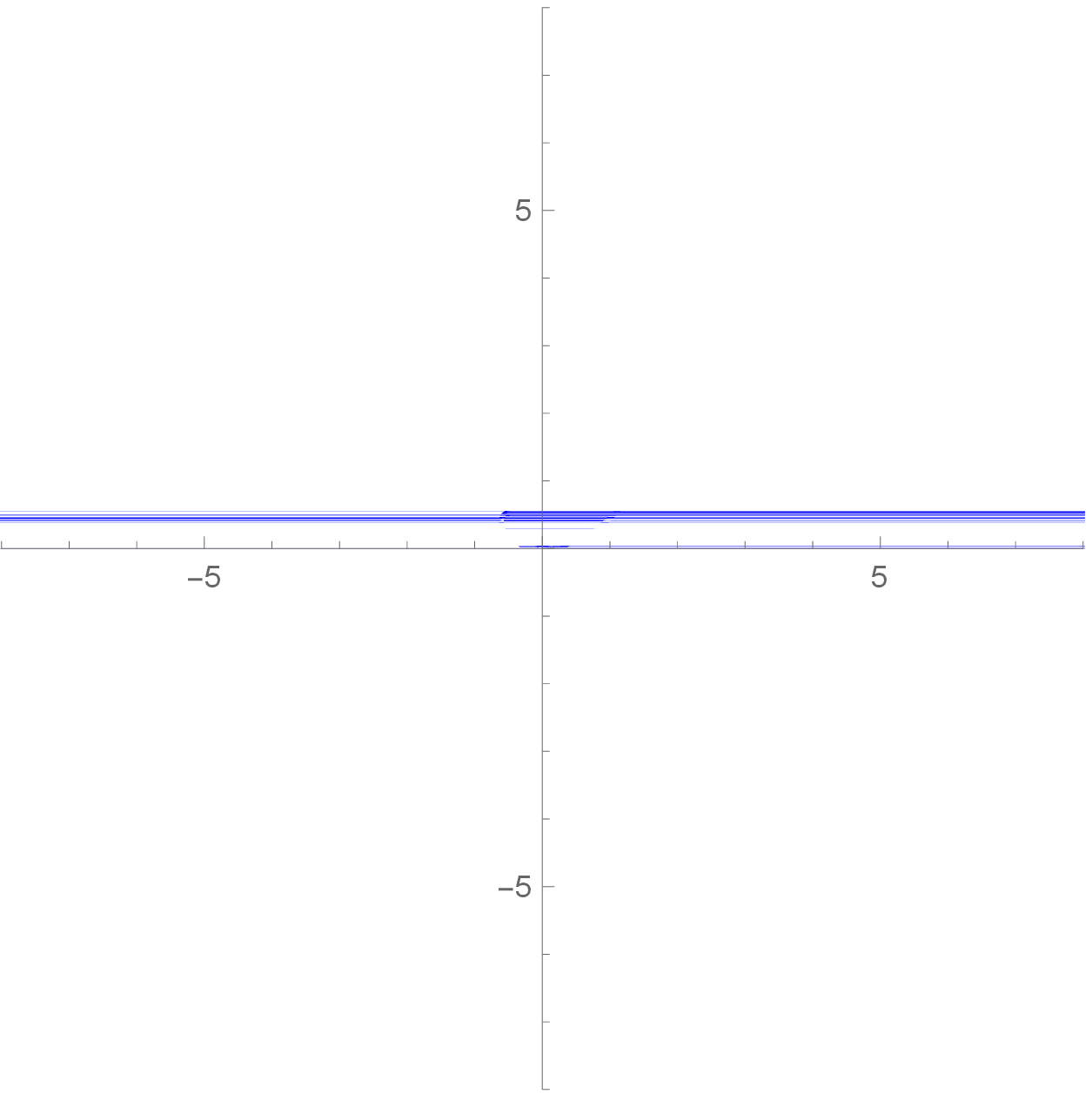}
		\caption{$\alpha=1$,\\Logarithmic spiral}\label{fig:log-aesthetic_alpha_1}
	\end{subfigure}
	\hfill
	\begin{subfigure}[t]{0.22\textwidth}
		\centering
		\includegraphics[width=\textwidth]{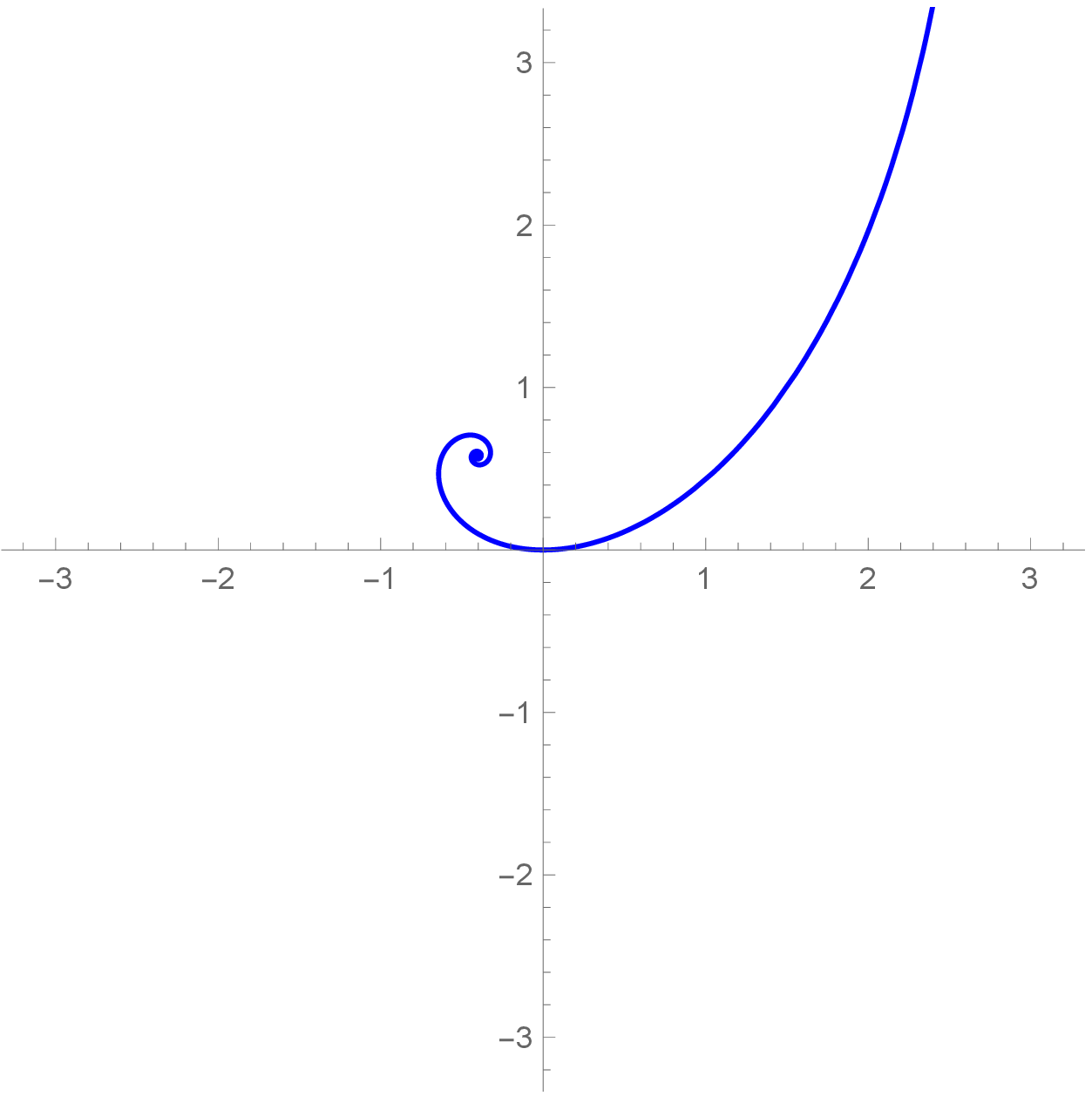}
		\caption{$\alpha=0.5$}\label{fig:log-aesthetic_alpha_0_5}
	\end{subfigure}
	\hfill 	
	\begin{subfigure}[t]{0.22\textwidth}
		\centering
		\includegraphics[width=\textwidth]{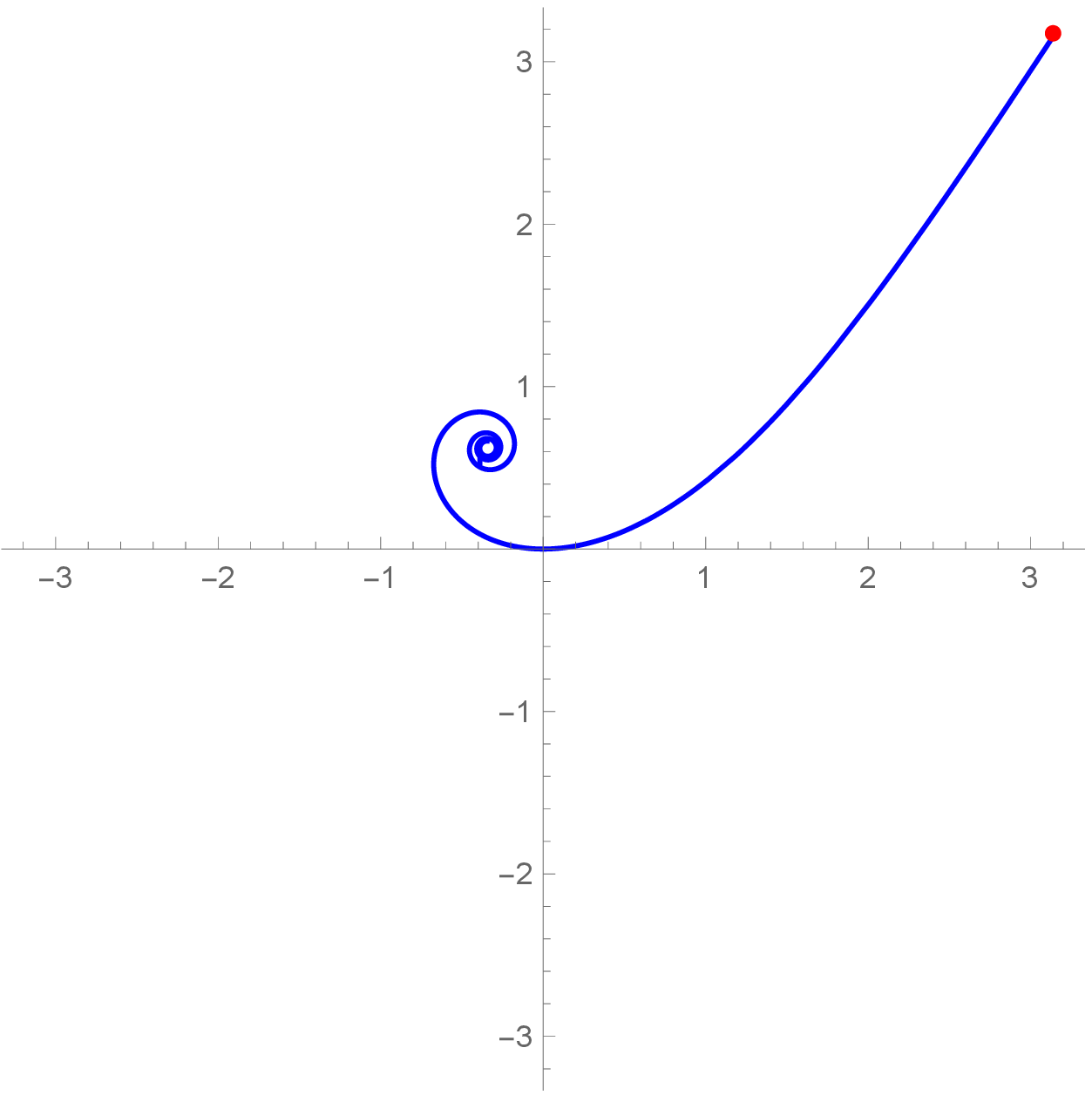}
		\caption{$\alpha=0$, \\ Nielsen's spiral}
	\end{subfigure}
	\hfill
	\begin{subfigure}[t]{0.22\textwidth}
		\centering
		\includegraphics[width=\textwidth]{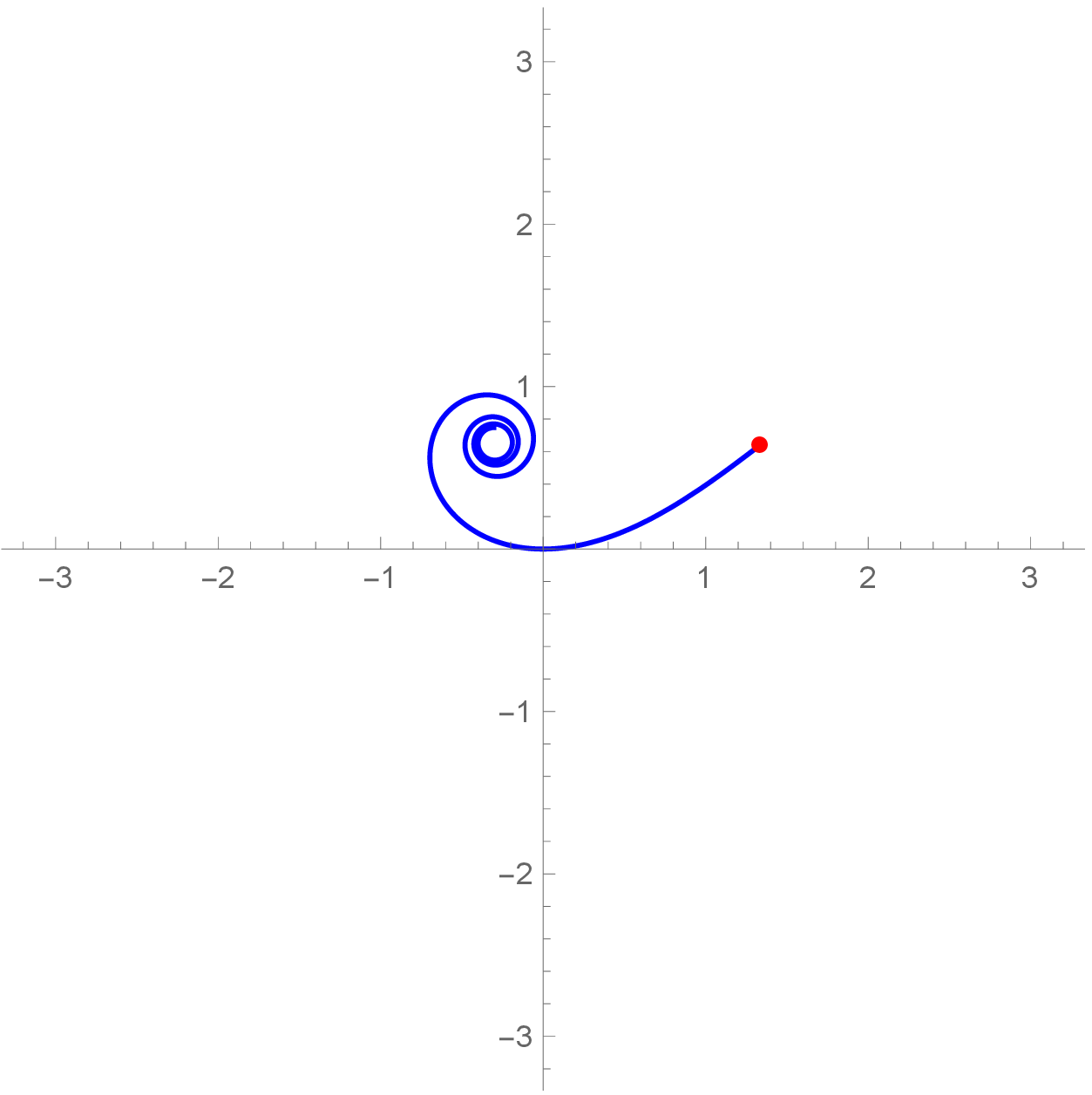}
		\caption{$\alpha=-0.5$}
	\end{subfigure}
	\\[2em]
	\centering
	\begin{subfigure}[t]{0.22\textwidth}
		\centering
		\includegraphics[width=\textwidth]{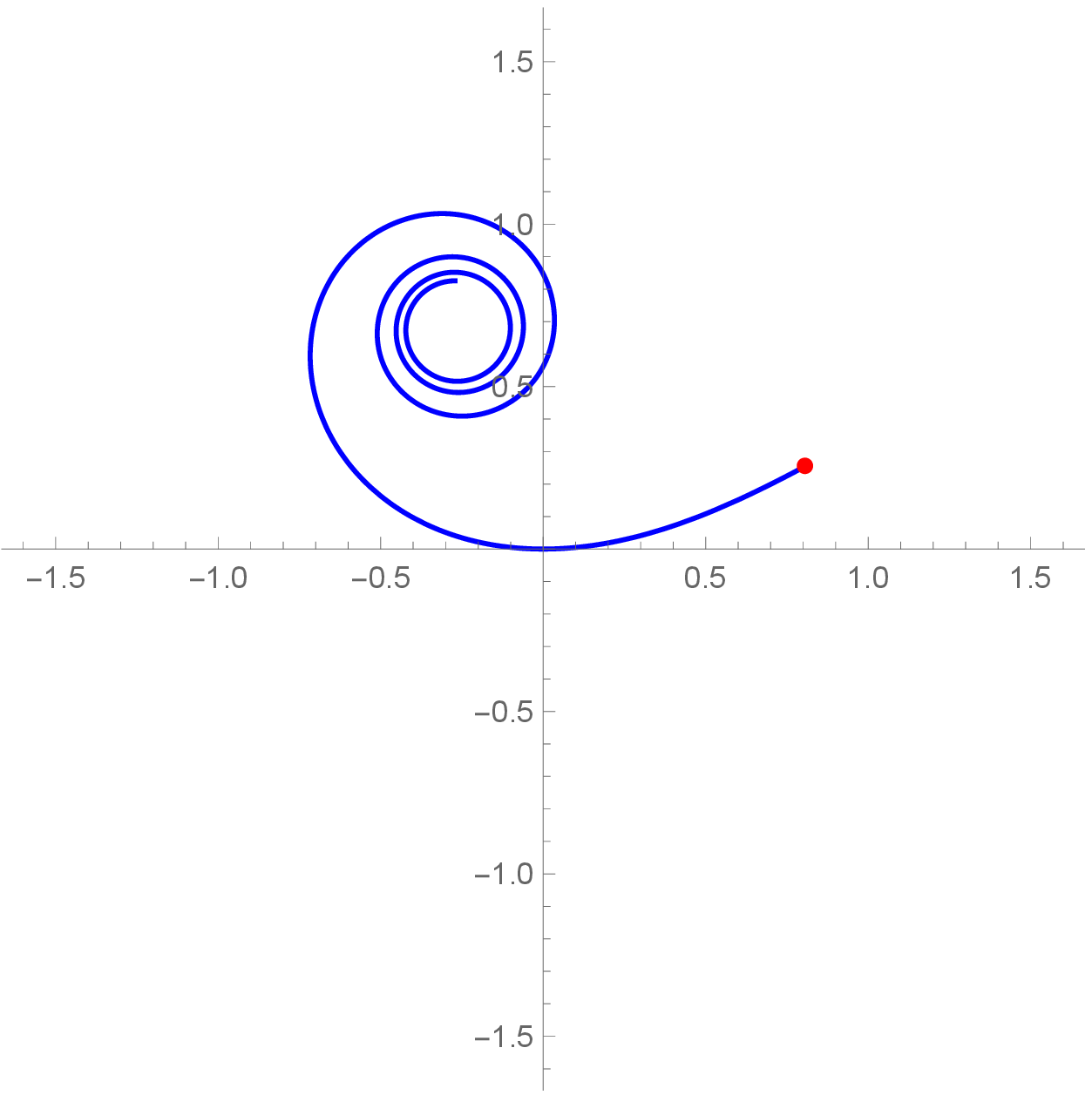}
		\caption{$\alpha=-1$,\\Clothoid}
	\end{subfigure}
 \hfill 	
	\begin{subfigure}[t]{0.22\textwidth}
		\centering
		\includegraphics[width=\textwidth]{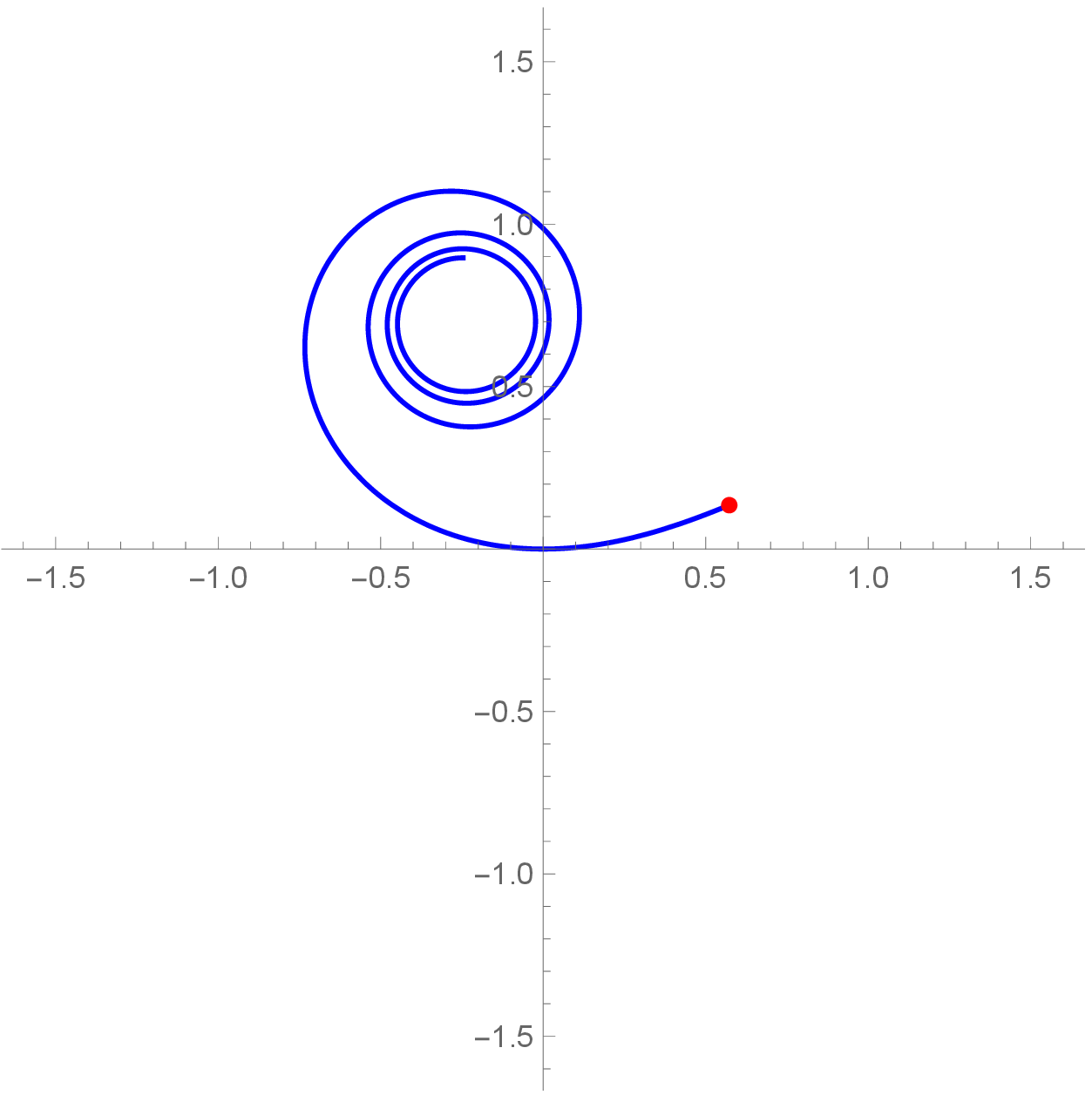}
		\caption{$\alpha=-1.5$}
	\end{subfigure}
 \hfill 	
	\begin{subfigure}[t]{0.22\textwidth}
		\centering
		\includegraphics[width=\textwidth]{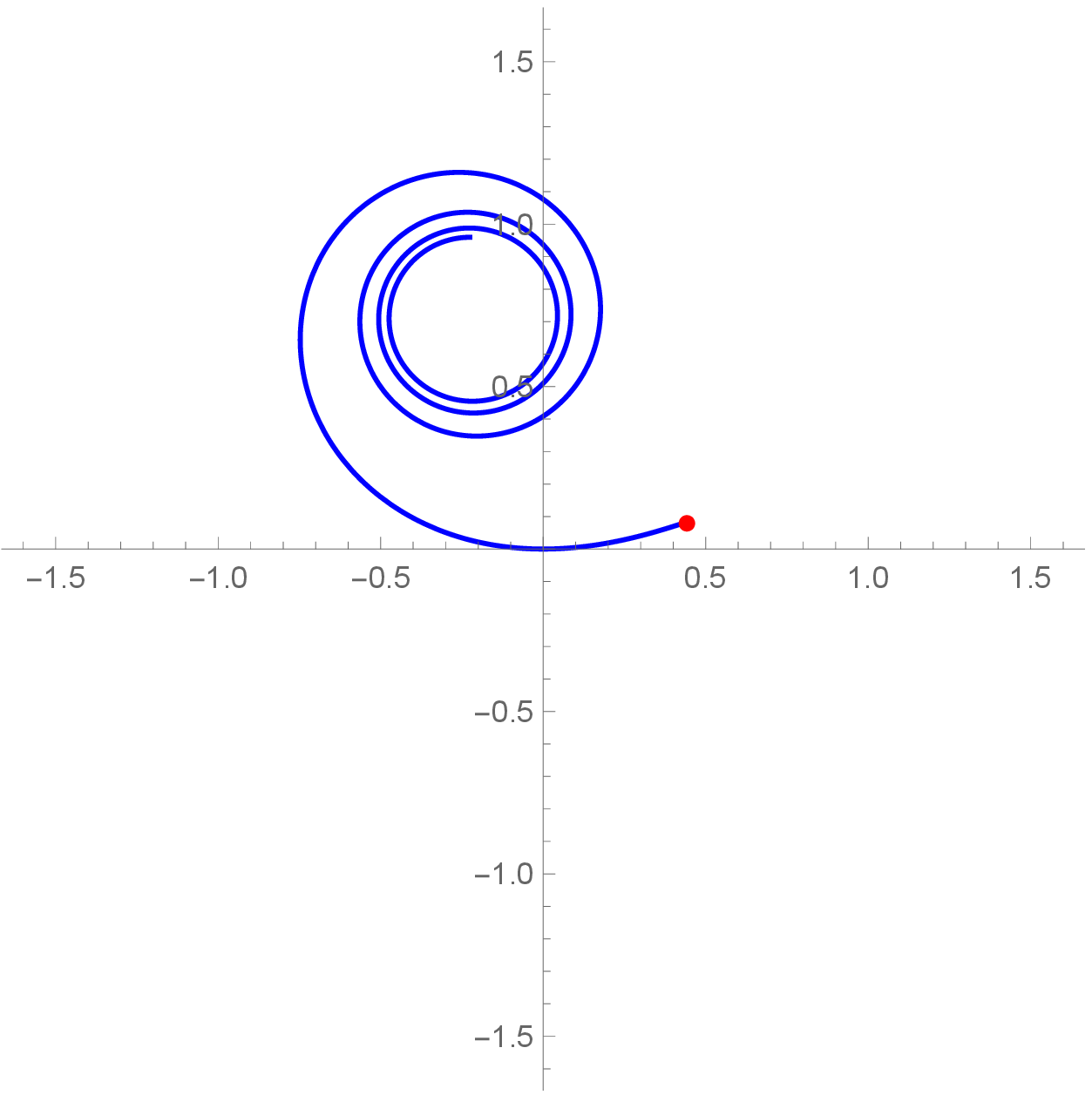}
		\caption{$\alpha=-2$}
	\end{subfigure}
 \hfill 	
	\begin{subfigure}[t]{0.22\textwidth}
		\centering
		\includegraphics[width=\textwidth]{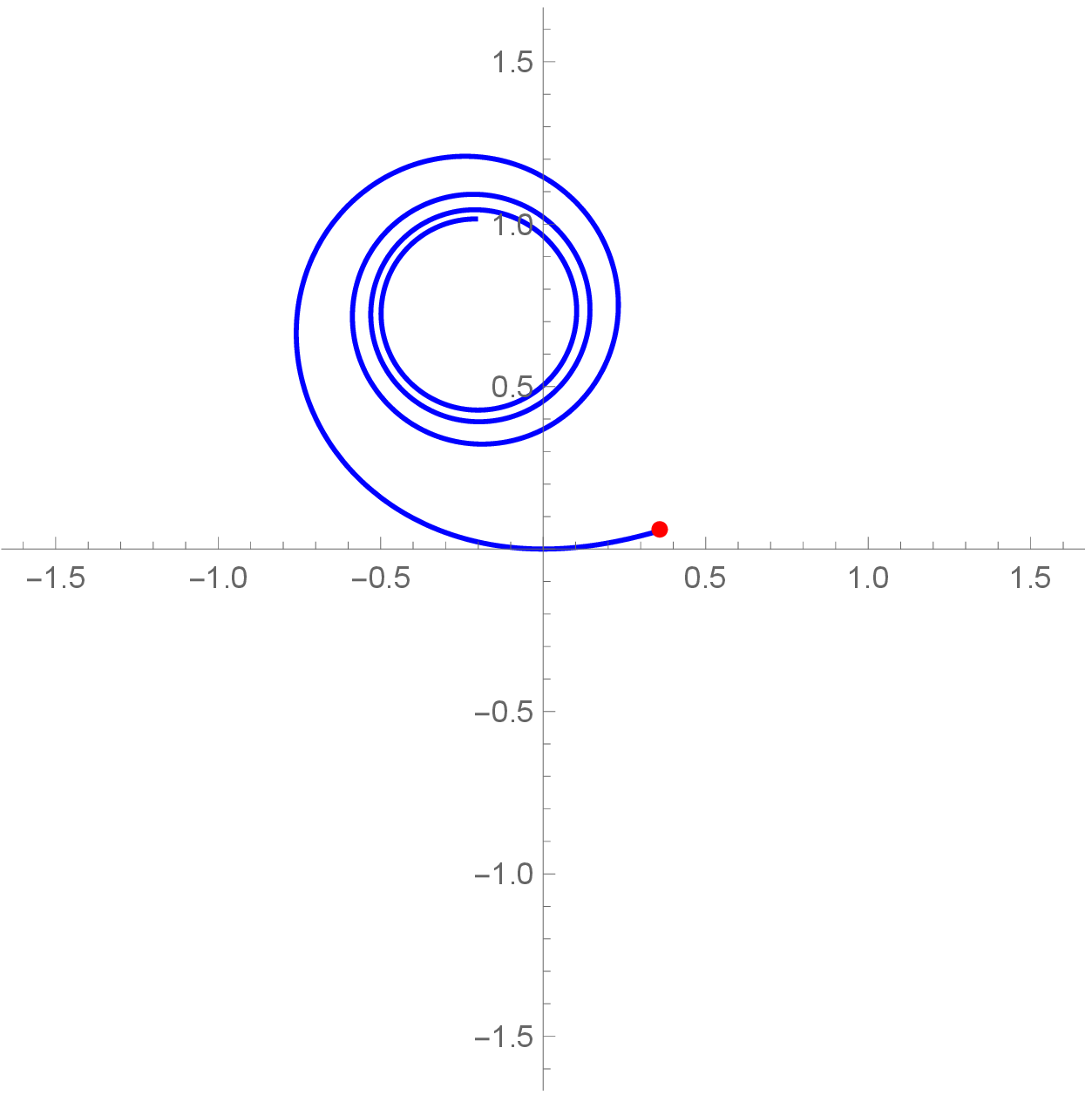}
		\caption{$\alpha=-2.5$}\label{fig:log-aesthetic_alpha_-2_5}
	\end{subfigure}
\caption{Log-aesthetic curves with different shape parameter $\alpha$ and $\Lambda=1$. The red points are the current bounds of $\theta$.}
\label{fig:log-aesthetic_alphas}
\end{figure}

Besides the general formula, the authors of \cite{yoshida2006interactive} also analyzed the properties of the curve. They called the case of $\Lambda=1$ as the standard form I of the curve. Let us highlight some important characteristics in this case regarding the isoptic of the curve, based on the value of $\alpha$:
\begin{description}
    \item[$\alpha>1$] As $\theta$ approaches $\infty$, the curve spirally diverges towards the point at $\rho=\infty$ (the arc length to that point is infinite). However, $\theta$ has a lower bound of $\frac{1}{\Lambda(1-\alpha)}$ (see \figref{fig:log-aesthetic_alpha_3}--\subref{fig:log-aesthetic_alpha_1_5}).
    \item[$\alpha<1$] As $\theta$ approaches $-\infty$, the curve spirally converges to the point at $\rho=0$ (the arc length to that point is infinite). However, $\theta$ has an upper bound of $\frac{1}{\Lambda(1-\alpha)}$ (see \figref{fig:log-aesthetic_alpha_0_5}--\subref{fig:log-aesthetic_alpha_-2_5}).
    \item[$\alpha=1$] As $\theta$ approaches $-\infty$, the curve spirally converges to the point at $\rho=0$ (the arc length to that point is finite), and as $\theta$ approaches $\infty$, the curve spirally diverges towards the point at $\rho=\infty$ (the arc length to that point is infinite, see \figref{fig:log-aesthetic_alpha_1}).
    \item[$\alpha=\rpm\infty$] The log-aesthetic curve is a circle centered at $\big[0$ $1\big]^T$.
\end{description}

\begin{table}[htb]
\centering
\begin{tabular}{|c c c c c c c c |}
\hline
    & \multicolumn{3}{c}{Tangential angle ($\theta$)} &
    & \multicolumn{3}{c|}{Arc length ($s$)} \\
\hline
    & $\alpha<1$ & $\alpha=1$ & $\alpha>1$ & &
      $\alpha<0$ & $\alpha=0$ & $\alpha>0$ \\
    Upper bound:  &
    $\frac{1}{\Lambda(1-\alpha)}$ & - & - & &
    $-\frac{1}{\Lambda\alpha}$ & - & - \\
    Lower bound:  &
    - & - & $\frac{1}{\Lambda(1-\alpha)}$ & &
    - & - & $-\frac{1}{\Lambda\alpha}$ \\
\hline
\end{tabular}
\caption[]{Upper and lower bound of $\theta$ and $s$.}
\label{table:arc_theta_bounds}
\end{table}

\subsection{Isoptic curves} \label{sec:isoptic_curve}

\begin{definition} \label{def:isoptic_curve}
For a given plane curve the isoptic curve is the locus of points from where the given curve can be seen under a predefined angle ${(0 < \gamma < \pi)} $\cite{yates1952}. The isoptic curve with right angle called orthoptic curve.
\end{definition}

\begin{figure}[h]
    \centering
    \includegraphics[width=0.44\textwidth]{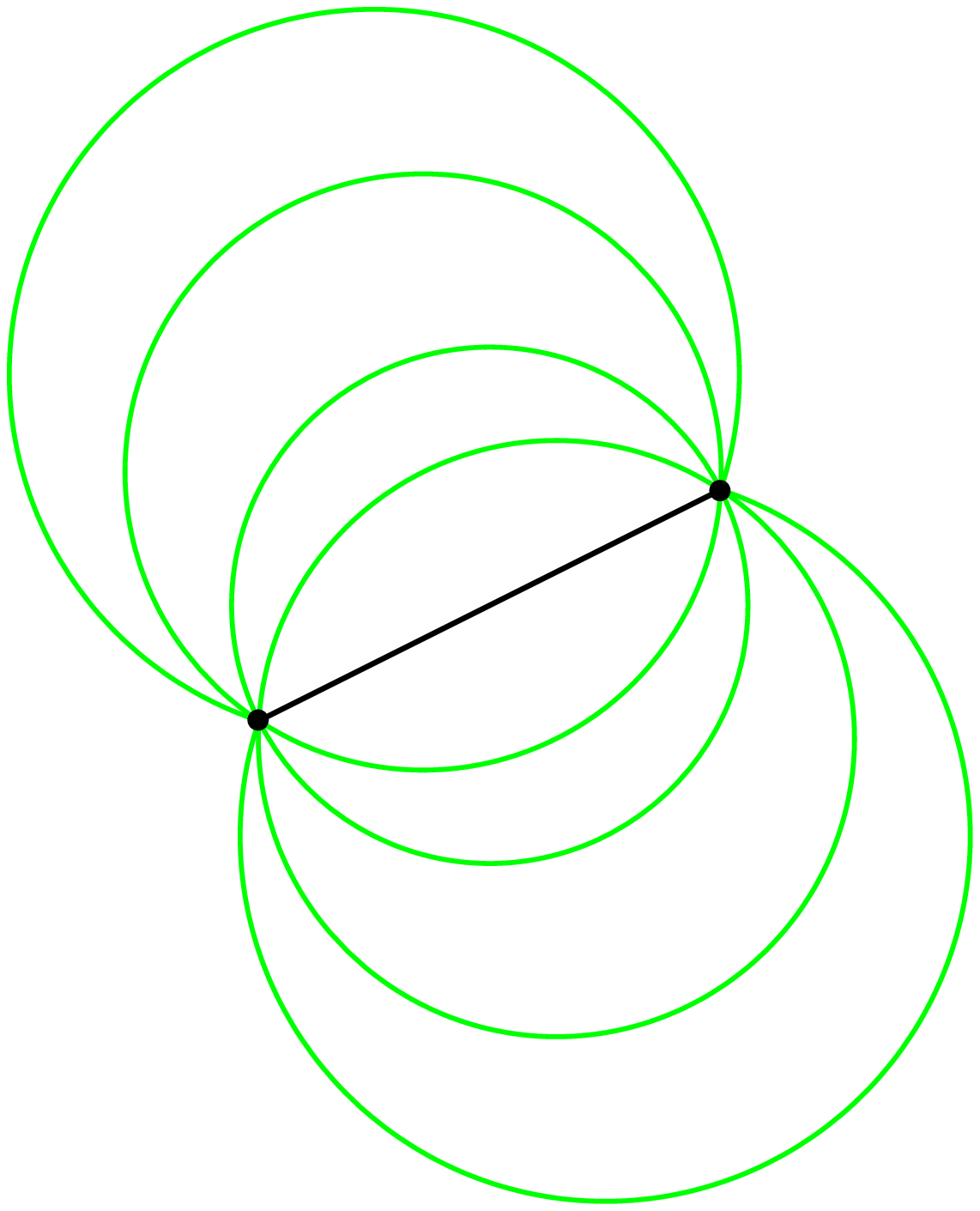}
    \hspace*{\fill}
    \includegraphics[width=0.44
    \textwidth]{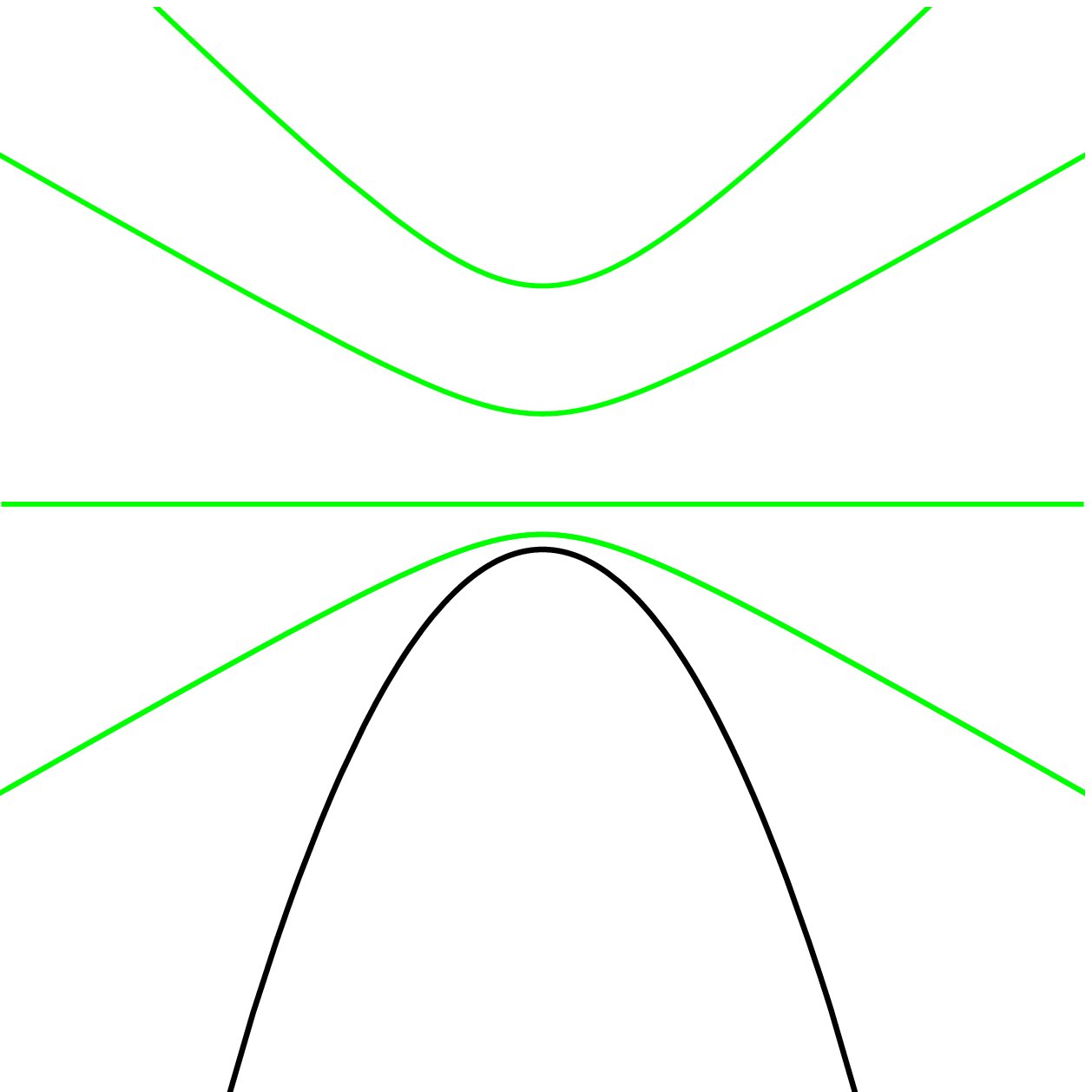}    
    \caption{Isoptic curves (green) of a line segment and a parabola with ${\gamma=\frac{3\pi}{4}} \textrm{, } {\frac{2\pi}{3}} \textrm{, } {\frac{\pi}{2}} \textrm{, and } {\frac{\pi}{3}}$. The orthoptic of the line segment is a circle, and the ortoptic of the parabola is its directrix.} \label{fig:classic-isoptic}
\end{figure}

For some planar, second-degree curves, the isoptic can be constructed simply based on the definition (see \figref{fig:classic-isoptic}, or for an overview in \cite{yates1952}). The authors of \cite{siebeck1860ueber} calculated isoptic curves of conic sections.
There are also interesting results of curves whose isoptic is a circle in \cite{wunderlich1971kurven,kurusa2012convex}, or whose isoptic is an ellipse in \cite{wunderlich1972kurven}.
The work of Cieślak et al. is also remarkable \cite{cieslak1991isoptics,cieslak1996isoptics}, in which, the authors define the isoptics of closed convex curves using support function.
Besides, there are also important results of isoptics in non-Euclidean geometry, e.g. in \cite{csima2013isoptic} or in \cite{csima2014isoptic}.

It should be noted that for higher degree curves the direct computation of the isoptic based on the definition can be extremely difficult. Therefore, Kunkli et al. defined a new calculation for free-form curves \cite{kunkli2013isoptics}. The authors determine the isoptic curve of convex Bézier-curves as the envelope of envelopes of families of isoptic circles over the chords of the curve.
Although, the calculations can easily be applied for other curve types, such as B-splines or NURBS, for complex or higher degree curves the isoptic usually can only be approximated.

\section{Isoptics of log-aesthetic curves} \label{sec:isoptic_log-aesthetic}

\begin{figure}[h]
  \centering
		\begin{subfigure}[t]{0.4\textwidth}
		\centering
		\includegraphics[width=\textwidth]{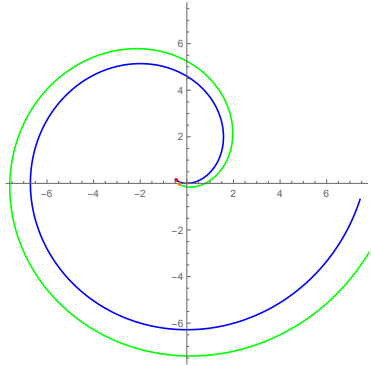}
		\caption{$\alpha=2$\\Circle involute}
	\end{subfigure}
 \hfill 	
	\begin{subfigure}[t]{0.4\textwidth}
		\centering
		\includegraphics[width=\textwidth]{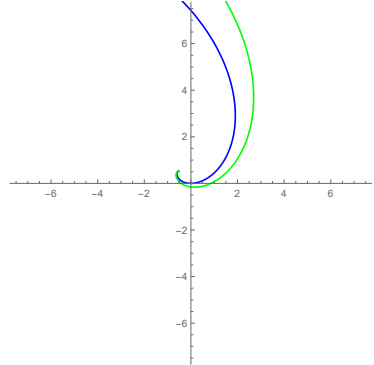}
		\caption{$\alpha=1$\\Logarithmic spiral}
	\end{subfigure}
	\\[1em]
	\begin{subfigure}[t]{0.4\textwidth}
		\centering
		\includegraphics[width=\textwidth]{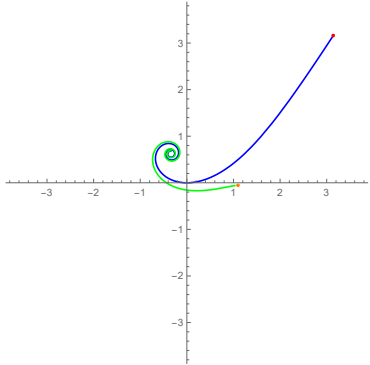}
		\caption{$\alpha=0$\\Nielsen’s spiral}
	\end{subfigure}
 \hfill 	
	\begin{subfigure}[t]{0.4\textwidth}
		\centering
		\includegraphics[width=\textwidth]{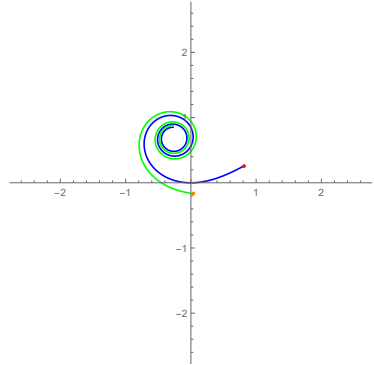}
		\caption{$\alpha=-1$\\Clothoid}
	\end{subfigure}
\caption{Log-aesthetic curves with $\Lambda=1$ (blue) and their isoptics with given $\delta=\frac{\pi}{3}$ (green), with the point of bounds (red and orange).}
\label{fig:log-aesthetic_alpha_isos}
\end{figure}

Based on \defref{def:isoptic_curve} and using the tangential angle parametrized \eqref{eq:log-aesthetic_point_theta} of the log-aesthetic curve the isoptic curve can be constructed as the  intersection of the appropriate involving lines. For an arbitrary point $P(\theta)$ of the log-aesthetic curve the corresponding curve point, which tangent line intersects the involving line of $P(\theta)$ under a given angle $\gamma$ is calculated as $P\left(\theta+(\pi-\gamma)\right)$. The tangent with unit length of the point $\theta$ is: $[\cos \theta, \sin \theta]^T$.
It is also true for all the tangential angle parametrized planar $r(t)$ curves \cite{abbena1998modern}:
\begin{align} \label{eq:log-aesthetic_unit_tangent}
     \frac{r^{\prime}(t)}{\left| r^{\prime}(t) \right|}
    	&= \begin{bmatrix}
    		\cos t \\
    		\sin t
    	\end{bmatrix}.
\end{align}
Let us define the isoptic with angle $\delta=\pi-\gamma$. Thus, the parametric equations of the appropriate tangent lines are
\begin{align} \begin{split}
		P(t_\theta) &= P(\theta) + t_\theta \begin{bmatrix}
           \cos \theta \\
           \sin \theta
         \end{bmatrix} \\
		P(t_{\theta+\delta}) &= P(\theta+\delta) + t_{\theta+\delta} \begin{bmatrix}
           \cos (\theta+\delta) \\
           \sin (\theta+\delta)
         \end{bmatrix}.
\end{split} \end{align}
Solving it for $t_\theta$, we obtain:
\begin{equation} \label{eq:log-aesthetic-isoptic_t}
	t_\theta = \csc (\delta) \left( \vec{V}^\delta_x(\theta) \sin(\theta+\delta) - \vec{V}^\delta_y(\theta) \cos(\theta+\delta) \right),
\end{equation}
where $\vec{V}^\delta_x(\theta)$ and $\vec{V}^\delta_y(\theta)$ are the coordinates of the vector 
$\vec{V}^\delta(\theta)=\overrightarrow{P(\theta)P(\theta+\delta)}$, which explicit equation can be written (from \eqref{eq:log-aesthetic_point_theta}) on the complex plane as
\begin{equation} \label{eq:log-aesthetic_vector_theta}
    \vec{V}^\delta(\theta)= 
\begin{cases}
    \int_{\theta}^{\theta+\delta}  \me^{(1+\iu)\Lambda\psi}\  d\psi
        & \text{if } \alpha = 1\\
    \int_{\theta}^{\theta+\delta}  \big((\alpha-1)\Lambda\psi+1\big)^{\frac{1}{\alpha-1}} \me^{\iu\psi}\  d\psi
        & \text{otherwise}.
\end{cases}
\end{equation}
Therefore, the isoptic curve of the log-aesthetic (and all the tangential angle parametrized) curve(s) with a given $\delta$ ($= \pi-\gamma$) angle can be calculated as (see \figref{fig:log-aesthetic_alpha_isos}):
\begin{equation} \label{eq:log-aesthetic_isoptic_curve}
I_{\delta}(\theta) = P(\theta)
	+ \csc (\delta) \left( \vec{V}^\delta_x(\theta) \sin(\theta+\delta) - \vec{V}^\delta_y(\theta) \cos(\theta+\delta) \right)
			\begin{bmatrix}
				 \cos \theta \\
				 \sin \theta
			 \end{bmatrix}
\end{equation}

The following modifications giving us more geometric meaning. The \eqref{eq:log-aesthetic-isoptic_t} can be altered using the harmonic addition theorem \cite{weisstein2021} as 
\begin{equation}
	t_\theta = \csc (\delta) \sign\left(\vec{V}^\delta_x(\theta)\right) \norm{\vec{V}^\delta(\theta)} \sin{\left(\theta+\delta+\arctan{\left(\frac{-\vec{V}^\delta_y(\theta)}{\vec{V}^\delta_x(\theta)}\right)}\right)}.
\label{eq:log-aesthetic-isoptic_t_haronic}
\end{equation}
Changing the $\arctan$ function to the well-known $\arctantwo$, that calculates the tangential angle of vector $\vec{V}^\delta(\theta)$ (a.k.a. the angle between the positive x axis and the ray to the point $\big[\vec{V}^\delta_x(\theta),\vec{V}^\delta_y(\theta)\big]^T$ and returns the value in the interval of $[0,\pi]$) the $\sign\left(\vec{V}^\delta_x(\theta)\right)$ can be omitted:
\begin{equation}
	t_\theta = \csc (\delta) \norm{\vec{V}^\delta(\theta)} \sin{\left(\theta+\delta-\arctantwo{\left(\vec{V}^\delta_y(\theta),\vec{V}^\delta_x(\theta)\right)}\right)}.
\label{eq:log-aesthetic-isoptic_t_simple}
\end{equation}
It means that the isoptic curve of the tangential angle parametrized curves can also be given as
\begin{equation} \label{eq:log-aesthetic_isoptic_curve_alt}
I_{\delta}(\theta) = P(\theta) + \csc (\delta) \norm{\vec{V}^\delta(\theta)} \sin{\left(\theta+\delta-\arctantwo{\left(\vec{V}^\delta_y(\theta),\vec{V}^\delta_x(\theta)\right)}\right)}
			\begin{bmatrix}
				 \cos \theta \\
				 \sin \theta
			 \end{bmatrix}
\end{equation}

To analyze the autoisoptic property, let us see the logarithmic curvature graph of the log-aesthetic curve, parametrized by $\theta$:
\begin{equation} \label{eq:LCG_theta}
P_{LCG}(\theta)=
	\begin{bmatrix}
		\log \rho(\theta) \\
		\log \frac{\rho(\theta)}{\mid\rho(s)\mid}
	\end{bmatrix},
\end{equation}
By substituting the appropriate $\rho$ functions, after some simplifications, it is:
\begin{equation} \label{eq:LCG_theta_exp}
P_{LCG}(\theta)=
\begin{cases}
	\begin{bmatrix}
		\Lambda \\
		\theta  \Lambda -\log (\Lambda )
	\end{bmatrix} & \alpha =1 \\[1em]
	\begin{bmatrix}
		\frac{\log ((\alpha -1) \theta  \Lambda +1)}{\alpha -1} \\
		\frac{\alpha  \log ((\alpha -1) \theta  \Lambda +1)}{\alpha -1}-\log (\Lambda )
	\end{bmatrix} & \text{otherwise.}
\end{cases}
\end{equation}
The \eqref{eq:LCG_theta_exp} on the LCG (which axes are $\log \rho$ and $\log \frac{\Delta s}{\Delta \rho / \rho}$) defines a straight line with slope $\alpha$ (see \figref{fig:LCG_all}).
The slope (which is constant for log-aesthetic curves) can be calculated as the following (by taking into account the \tableref{table:arc_theta_bounds} limits of $\theta$):
\begin{equation} \label{eq:LCG_slope}
\hat{\alpha}(\theta) =
\begin{cases}
	\frac{P_{LCG}(\theta-\phi)_y-P_{LCG}(\theta)_y}{P_{LCG}(\theta-\phi)_x-P_{LCG}(\theta)_x}
		& \alpha < 1 \\[1em]
	\frac{P_{LCG}(\theta+\phi)_y-P_{LCG}(\theta)_y}{P_{LCG}(\theta+\phi)_x-P_{LCG}(\theta)_x}
		& \alpha \geq 1,
\end{cases}
\end{equation}
where $\phi \in \{\mathbb{R}\setminus{0}\}$ and $P_{LCG}(\theta)_x$ and $P_{LCG}(\theta)_y$ is the coordinates of the point on the LCG defined by $\theta$.

Unfortunately, the exact logaritmic curvature graph of the isoptic curve, using the $\theta$ parametrized \eqref{eq:log-aesthetic_isoptic_curve} or \eqref{eq:log-aesthetic_isoptic_curve_alt} cannot be determined for arbitrary $\alpha \in \ \mathbb{R}$ (not even using computer algebra softwares). Solely, in the cases when the log-aesthetic curve is the logarithmic spiral ($\alpha=1$) or when the \eqref{eq:log-aesthetic_point_theta} can be represented in terms of trigonometric functions ($\alpha=2, \frac{3}{2}, \frac{4}{3}, \dots , \frac{k+1}{k}$, $k \in \mathbb{N}^{*}$) or using Fresnel integrals ($\alpha=-1, 3$) \cite{ziatdinov2012analytic}.

\begin{figure}[h]
\centering
  \includegraphics[width=0.9\linewidth]{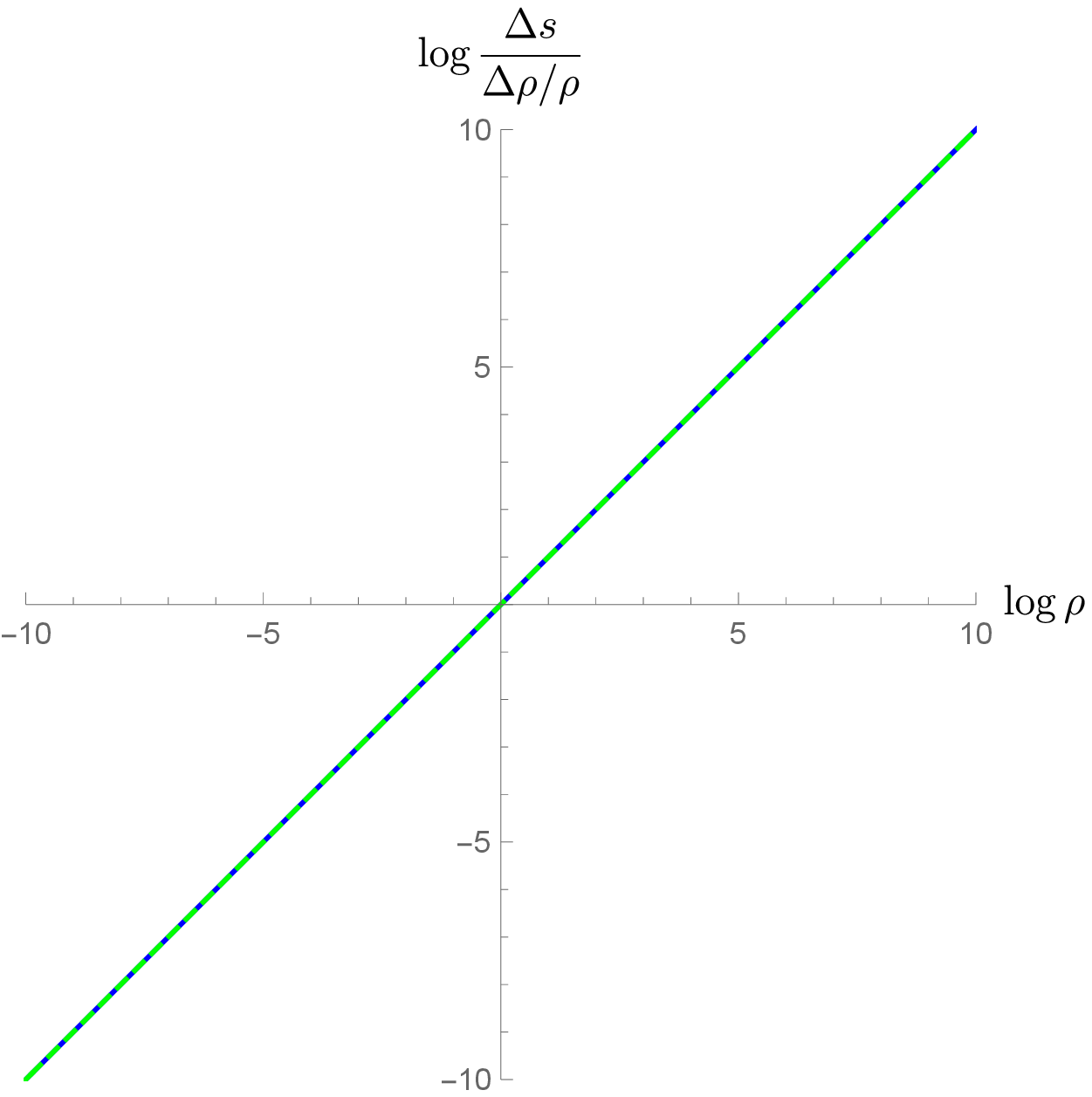}
\caption{Logarithmic curvature graph of the log-aesthetic curve with $\alpha=1$, $\Lambda=1$ (blue) and its isoptic with $\delta=\frac{2}{3}\pi$ (green, dashed).}
\label{fig:LCG_iso_1}
\end{figure}

Let us see the case of $\alpha=1$. The LCG equation of the log-aestheic curve is:
\begin{equation} \label{eq:LCG_1}
P^{\alpha=1}_{LCG}(\theta)=
	\begin{bmatrix}
		\theta  \Lambda \\
		\theta  \Lambda -\log (\Lambda )
	\end{bmatrix},
\end{equation}
and the LCG formula of its isoptic is
\begin{equation} \label{eq:LCG_1_iso}
I^{\alpha=1}_{LCG}(\theta)=
	\begin{bmatrix}
		\frac{1}{2} \log \left(\frac{e^{2 \delta  \Lambda }-2 e^{\delta  \Lambda } \cos (\delta )+1}{\Lambda ^2+1}\right)+\log (\csc (\delta ))+\theta  \Lambda \\
		\log \left(\frac{\csc (\delta ) \sqrt{\frac{e^{2 \delta  \Lambda }-2 e^{\delta  \Lambda } \cos (\delta )+1}{\Lambda ^2+1}}}{\Lambda }\right)+\theta  \Lambda
	\end{bmatrix}.
\end{equation}
By calculating the slope $\hat{\alpha}(\theta)$ of the isoptic, using a similar formula as \eqref{eq:LCG_slope}, we see that it is constant $1$ for all $\theta \in \R$, meaning that, it is also a log-aesthetic curve with shape parameter 1. It proves its autoisoptic property (see \figref{fig:LCG_iso_1}).

For another example, let us see the LCG of the log-aesthetic curve with $\alpha=2$ (circle involute), that is
\begin{equation} \label{eq:LCG_2}
P^{\alpha=2}_{LCG}(\theta)=
	\begin{bmatrix}
		\log (\theta  \Lambda +1) \\
		2 \log (\theta  \Lambda +1)-\log (\Lambda )
	\end{bmatrix},
\end{equation}
and the LCG of the isoptic with $\delta=\frac{2\pi}{3}$ and $\Lambda=1$ is (note: it is not the general formula with arbitrary $\Lambda>0$ and ${0 < \delta < \pi}$ because it is too long to appear in print):
\begin{equation} \label{eq:LCG_2_iso}
I^{\alpha=2,\Lambda=1}_{LCG,\delta=2\pi/3}(\theta)=
	\begin{bmatrix}
		\log \left(\frac{2 \left(-6 \pi  \left(-3 \theta +\sqrt{3}-3\right)+27 \theta  (\theta +2)+4 \pi ^2+54\right)^{3/2}}{27 \sqrt{3} \left| 3 \theta  (\theta +2)+\pi  \left(2 \theta -\sqrt{3}+2\right)+\frac{4 \pi ^2}{9}+9\right| }\right)\\[2.5em]
		\log \left(\frac{2 \left(-6 \pi  \left(-3 \theta +\sqrt{3}-3\right)+27 \theta  (\theta +2)+4 \pi ^2+54\right)^{3/2}}{27 \sqrt{3}}\right.+\\
		+\left.\frac{\left| \frac{-9 \pi  \left(-2 \theta +\sqrt{3}-2\right)+27 \theta  (\theta +2)+4 \pi ^2+81}{(3 \theta +\pi +3) \left(3 \pi  \left(6 \theta -5 \sqrt{3}+6\right)+27 (\theta  (\theta +2)+5)+4 \pi ^2\right)}\right| }{27
   \sqrt{3}}\right)
	\end{bmatrix}.
\end{equation}

\begin{figure}[htb]
\centering
  \includegraphics[width=0.9\linewidth]{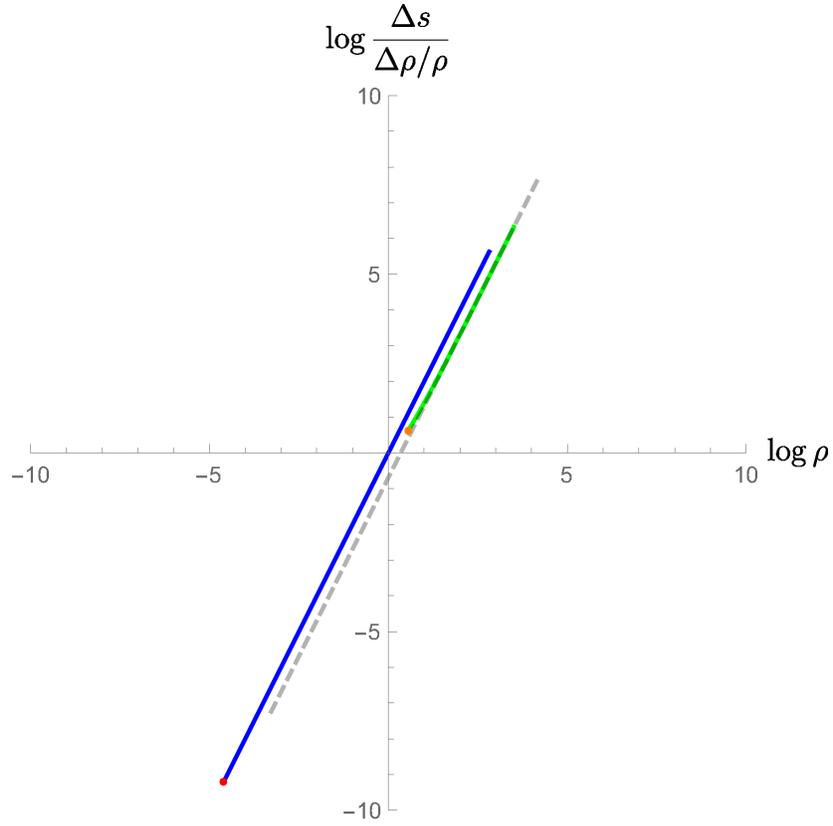}
\caption{Logarithmic curvature graph of the log-aesthetic curve with $\alpha=2$, $\Lambda=1$ (blue and the dashed, that is a transformed) and its isoptic with $\delta=\frac{2}{3}\pi$ (green).}
\label{fig:LCG_iso_2}
\end{figure}

By calculating its $\hat{\alpha}(\theta)$ slope, similarly as \eqref{eq:LCG_slope} with different $\theta$ values (and with $\phi=\pi$), we see: 
\begin{equation} \label{eq:LCG_iso_alpha_2}
\begin{split}
    \hat{\alpha}(-1) &\approx 1.86204 \text{ (lower bound of } \theta  \text{)}\\
    \hat{\alpha}(0) &\approx 1.96486\\
    \hat{\alpha}(\pi) &\approx 1.99523\\
    \hat{\alpha}(3\pi) &\approx 1.99899\\
    \hat{\alpha}(10\pi) &\approx 1.99987\\
    \hat{\alpha}(100\pi) &\approx 1.99999845\\
    \lim_{\theta\to\infty} \hat{\alpha}(\theta) &= 2.
\end{split}
\end{equation}
Therefore, it can be claimed that the log-aesthetic curve with $\alpha=2$ does not coincide with its isoptic, meaning that the circle involute is not autoisoptic (see \figref{fig:LCG_iso_2}). However, the limit of $\hat{\alpha}(\theta)$ is 2, when $\theta$ approaches $\infty$.

We receive similar results by calculating and drawing the different logarithmic curvature graphs of the isoptic curve in the other mentioned $\alpha>1$ cases. The slope of the isoptic is not constant, but the limit of the slope as the tangential angle approaches infinity is always equal to the specified $\alpha$ parameter. 
We can draw the following conclusion: the log-aesthetic curve is not autoisoptic in all the cases of $\alpha>1$.

\begin{figure}[htb]
\centering
  \includegraphics[width=0.9\linewidth]{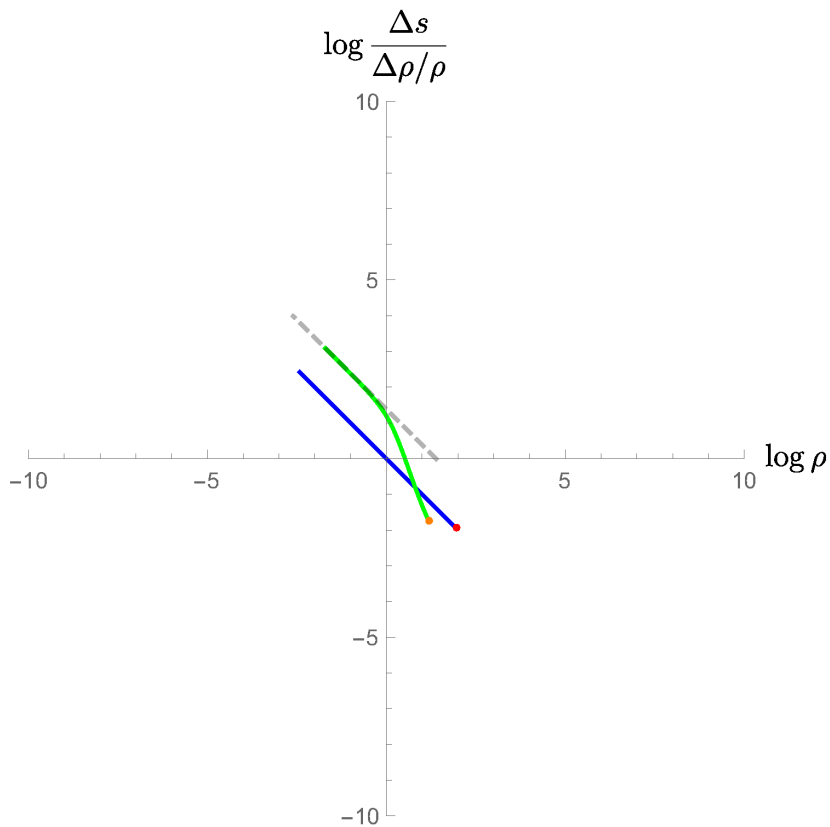}
\caption{Logarithmic curvature graph of the log-aesthetic curve with $\alpha=-1$, $\Lambda=1$ (blue and the dashed, that is a transformed) and its isoptic with $\delta=\frac{2}{3}\pi$ (green).}
\label{fig:LCG_iso_-1}
\end{figure}

Let us see another example, the log-aesthetic curve with $\alpha=-1$ (clothoid). In this case,
the general equation of the log-aesthetic curve is defined using Fresnel integrals \cite{ziatdinov2012analytic}.
The LCG of it can be obtained as
\begin{equation} \label{eq:LCG_-1}
P^{\alpha=-1}_{LCG}(\theta)=
	\begin{bmatrix}
		-\frac{1}{2} \log (1-2 \theta  \Lambda ) \\
		\frac{1}{2} \log (1-2 \theta  \Lambda )-\log (\Lambda )
	\end{bmatrix}.
\end{equation}
However, the LCG of its isoptic can be obtained with a fixed $\Lambda$ and $\delta$ values using computer algebra systems, it is too long to appear in print.
We calculated it for $\Lambda=1$ and $\delta=\frac{2}{3}\pi$ and computed the $\hat{\alpha}(\theta)$ slope of it with different $\theta$ values and with $\phi=\pi$, similarly as in \eqref{eq:LCG_slope}. From the following cases
\begin{equation} \label{eq:LCG_iso_alpha_-1}
\begin{split}
    \hat{\alpha}(\frac{1}{2}-\delta) &\approx -2.23102 \text{ (upper bound of } \theta  \text{ for the isoptic)}\\
    \hat{\alpha}(-\pi) &\approx -1.27031\\
    \hat{\alpha}(-2\pi) &\approx -1.0705\\
    \hat{\alpha}(-5\pi) &\approx -1.01251\\
    \hat{\alpha}(-10\pi) &\approx -1.00329\\
    \hat{\alpha}(-100\pi) &\approx -1.00003\\
    \lim_{\theta\to-\infty} \hat{\alpha}(\theta) &= -1,
\end{split}
\end{equation}

we see that the log-aesthetic curve with $\alpha=-1$ does not coincide with its isoptic, meaning that the clothoid is not autoisoptic (see \figref{fig:LCG_iso_-1}). However, the limit of $\hat{\alpha}(\theta)$ is -1, when $\theta$ approaches $-\infty$.
Note, the slope $\hat{\alpha}(\theta)$ of the LGC of the isoptic in this case differs more from $\alpha$ around the upper bound rather than in the case of $\alpha=2$. The reason can be that fact the log aesthetic curve has an inflection point at the upper bound of $\theta$ in case of $\alpha<1$.

\section{Conclusion}

We have defined the isoptics of the log-aesthetic curve in the two-dimensional plane. We have also given an exact formula to obtain the isoptics of all the tangential angle parametrized plane curves.

We have examined three instances of the log-aesthetic curve based on its shape parameter. In the first case, when it is the logarithmic spiral ($\alpha=1$), we have shown its autoisoptic property.
In the second case, when it is the circle involute ($\alpha=2$) and in the third case, when it is the clothoid ($\alpha=-1$), we have proven that they are not autoisoptic curves. It means that the log-aesthetic curve with arbitrary $\alpha \in \mathbb{R}$ is not autoisoptic.

Moreover, by analyzing \eqref{eq:log-aesthetic-isoptic_t_simple}, we believe that the log-aesthetic curve coincides with its isoptic only in the cases of $\alpha=1$ and $\alpha=\rpm\infty$ because in this cases the tangential angle $\theta$ and the tangential angle of the vector $\vec{V}^\delta(\theta)$ are in a relation of proportionality ($\theta \propto \arctantwo{\left(\vec{V}^\delta_y(\theta),\vec{V}^\delta_x(\theta)\right)}$), such that the $\sin{\left(\theta+\delta-\arctantwo{\left(\vec{V}^\delta_y(\theta),\vec{V}^\delta_x(\theta)\right)}\right)}$ is constant and the $\norm{\vec{V}^\delta(\theta)}$ equals to the radius of curvature of the log-aesthetic curve with a constant multiplier. In the other cases, when $-\infty<\alpha<1$ or $1<\alpha<\infty$ the ratio of the tangential angle $\theta$ and the tangential angle of the vector $\vec{V}^\delta(\theta)$ is not constant, therefore the log-aesthetic curve does not fulfill these properties. We assume that this is the reason why the log-aesthetic curve is not autoisoptic in these cases.

\bibliographystyle{plain}
\bibliography{main}

\begin{thebibliography}{10}

\bibitem{abbena1998modern}
Elsa Abbena, Simon Salamon, and Alfred Gray.
\newblock {\em Modern differential geometry of curves and surfaces with
  Mathematica}.
\newblock CRC press, 1998.

\bibitem{cieslak1991isoptics}
Waldemar Cie{\'{s}}lak, Andrzej Miernowski, and Witold Mozgawa.
\newblock Isoptics of a closed strictly convex curve.
\newblock In {\em Global Differential Geometry and Global Analysis}, pages
  28--35. Springer Berlin Heidelberg, 1991.

\bibitem{cieslak1996isoptics}
Waldemar Cie{\'{s}}lak, Andrzej Miernowski, and Witold Mozgawa.
\newblock Isoptics of a closed strictly convex curve.-ii.
\newblock {\em Rendiconti del Seminario Matematico della Universit{\`a} di
  Padova}, 96:37--49, 1996.

\bibitem{csima2014isoptic}
G{\'e}za Csima and Szirmai Jen{\H{o}}.
\newblock Isoptic curves of conic sections in constant curvature geometries.
\newblock {\em Mathematical Communications}, 19(2):277--290, 2014.

\bibitem{csima2013isoptic}
G{\'e}za Csima and Jen{\H{o}} Szirmai.
\newblock On the isoptic hypersurfaces in the n-dimensional euclidean space.
\newblock {\em KoG (Scientific and professional journal of Croatian Society for
  Geometry and Graphics)}, 17(17.):53--57, 2013.

\bibitem{harada1999aesthetic}
T.~{Harada}, F.~{Yoshimoto}, and M.~{Moriyama}.
\newblock An aesthetic curve in the field of industrial design.
\newblock In {\em Proceedings 1999 IEEE Symposium on Visual Languages}, pages
  38--47. IEEE, 1999.

\bibitem{harada1997study}
Toshinobu Harada.
\newblock Study of quantitative analysis of the characteristics of a curve.
\newblock {\em Forma}, 12(1):55--63, 1997.

\bibitem{kunkli2013isoptics}
Roland Kunkli, Ildik{\'o} Papp, and Mikl{\'o}s Hoffmann.
\newblock Isoptics of b{\'e}zier curves.
\newblock {\em Computer Aided Geometric Design}, 30(1):78--84, 2013.

\bibitem{kurusa2012convex}
{\'A}rp{\'a}d Kurusa.
\newblock Is a convex plane body determined by an isoptic?
\newblock {\em Beitr{\"a}ge zur Algebra und Geometrie/Contributions to Algebra
  and Geometry}, 53(1):281--294, 2012.

\bibitem{levien2009interpolating}
Raph Levien and Carlo~H. S{\'e}quin.
\newblock Interpolating splines: which is the fairest of them all?
\newblock {\em Computer-Aided Design and Applications}, 6(1):91--102, 2009.

\bibitem{miura2006general}
Kenjiro~T. Miura.
\newblock A general equation of aesthetic curves and its self-affinity.
\newblock {\em Computer-Aided Design and Applications}, 3(1-4):457--464, 2006.

\bibitem{miura2005derivation}
Kenjiro~T Miura, Junji Sone, Atsushi Yamashita, and Toru Kaneko.
\newblock Derivation of a general formula of aesthetic curves.
\newblock In {\em Proceedings of the Eighth International Conference on Humans
  and Computers (HC2005)}, pages 166--171, 2005.

\bibitem{boris2019examples}
Boris Odehnal.
\newblock Examples of autoisoptic curves.
\newblock In Luigi Cocchiarella, editor, {\em ICGG 2018 - Proceedings of the
  18th International Conference on Geometry and Graphics}, pages 350--361,
  Cham, 2019. Springer International Publishing.

\bibitem{siebeck1860ueber}
Paul Siebeck.
\newblock Ueber eine gattung von curven vierten grades, welche mit den
  elliptischen functionen zusammenh{\"a}ngen.
\newblock {\em Journal für die reine und angewandte Mathematik},
  1860(57):359--370, 1860.

\bibitem{weisstein2021}
Eric~W. Weisstein.
\newblock Harmonic addition theorem. {From MathWorld---A Wolfram Web Resource}.
\newblock Visited: 12/01/2021.

\bibitem{wunderlich1971kurven}
Walter Wunderlich.
\newblock Kurven mit isoptischem kreis.
\newblock {\em Aequationes mathematicae}, 6(1):71--81, 1971.

\bibitem{wunderlich1972kurven}
Walter Wunderlich.
\newblock Kurven mit isoptischer ellipse.
\newblock {\em Monatshefte f{\"u}r Mathematik}, 75(4):346--362, 1972.

\bibitem{yates1952}
R.~C. Yates.
\newblock {\em {A Handbook on Curves and their Properties}}.
\newblock Ann Arbor, J.W. Edwards, 1947.

\bibitem{yoshida2009log}
Norimasa Yoshida, Ryo Fukuda, and Takafumi Saito.
\newblock Log-aesthetic space curve segments.
\newblock In {\em 2009 SIAM/ACM Joint Conference on Geometric and Physical
  Modeling}, pages 35--46, 2009.

\bibitem{yoshida2006interactive}
Norimasa Yoshida and Takafumi Saito.
\newblock Interactive aesthetic curve segments.
\newblock {\em The Visual Computer}, 22(9-11):896--905, 2006.

\bibitem{yoshida2012evolutes}
Norimasa Yoshida and Takafumi Saito.
\newblock The evolutes of log-aesthetic planar curves and the drawable
  boundaries of the curve segments.
\newblock {\em Computer-Aided Design and Applications}, 9(5):721--731, 2012.

\bibitem{ziatdinov2012analytic}
Rushan Ziatdinov, Norimasa Yoshida, and Tae wan Kim.
\newblock Analytic parametric equations of log-aesthetic curves in terms of
  incomplete gamma functions.
\newblock {\em Computer Aided Geometric Design}, 29(2):129--140, 2012.

\end{thebibliography}

\end{document}